\documentclass{article}
\usepackage{amsmath}
\usepackage{amssymb}
\begin{document}
\def\e#1\e{\begin{equation}#1\end{equation}}
\def\ea#1\ea{\begin{align}#1\end{align}}
\def\eq#1{{\rm(\ref{#1})}}
\newtheorem{thm}{Theorem}[section]
\newtheorem{prop}[thm]{Proposition}
\newtheorem{cor}[thm]{Corollary}
\newenvironment{dfn}{\medskip\refstepcounter{thm}
\noindent{\bf Definition \thesection.\arabic{thm}\ }}{\medskip}
\newenvironment{proof}[1][,]{\medskip\ifcat,#1
\noindent{\it Proof.\ }\else\noindent{\it Proof of #1.\ }\fi}
{\hfill$\square$\medskip}
\def\dim{\mathop{\rm dim}}
\def\Re{\mathop{\rm Re}}
\def\Im{\mathop{\rm Im}}
\def\vol{\mathop{\rm vol}}
\def\Tr{\mathop{\rm Tr}}
\def\U{\mathbin{\rm U}}
\def\SU{\mathop{\rm SU}}
\def\sn{{\textstyle\mathop{\rm sn}}}
\def\cn{{\textstyle\mathop{\rm cn}}}
\def\dn{{\textstyle\mathop{\rm dn}}}
\def\sech{{\textstyle\mathop{\rm sech}}}
\def\ge{\geqslant} 
\def\le{\leqslant} 
\def\cal{\mathcal}
\def\R{\mathbin{\mathbb R}}
\def\Z{\mathbin{\mathbb Z}}
\def\C{\mathbin{\mathbb C}}
\def\CP{\mathbb{CP}}
\def\al{\alpha}
\def\be{\beta}
\def\ga{\gamma}
\def\de{\delta}
\def\ep{\epsilon}
\def\th{\theta}
\def\la{\lambda}
\def\vp{\varphi}
\def\si{\sigma}
\def\De{\Delta}
\def\La{\Lambda}
\def\Om{\Omega}
\def\Si{\Sigma}
\def\om{\omega}
\def\d{{\rm d}}
\def\pd{\partial}
\def\ts{\textstyle}
\def\w{\wedge}
\def\lt{\ltimes}
\def\sm{\setminus}
\def\ov{\overline}
\def\ot{\otimes}
\def\iy{\infty}
\def\ra{\rightarrow}
\def\t{\times}
\def\op{\oplus}
\def\ms#1{\vert#1\vert^2}
\def\md#1{\vert #1 \vert}
\def\bmd#1{\bigl\vert #1 \bigr\vert}
\def\an#1{\langle#1\rangle}
\def\ban#1{\bigl\langle#1\bigr\rangle}
\title{Ruled special Lagrangian 3-folds in $\C^3$}
\author{Dominic Joyce, \\ Lincoln College, Oxford}
\date{December 2000}
\maketitle

\section{Introduction}
\label{ru1}

This is the fourth in a series of papers \cite{Joyc1,Joyc2,Joyc3}
constructing explicit examples of special Lagrangian submanifolds 
(SL $m$-folds) in $\C^m$. This paper focusses on {\it ruled}\/ special
Lagrangian 3-folds in $\C^3$, that is, SL 3-folds $N$ in $\C^3$ admitting
a fibration $\pi:N\ra\Si$ for some 2-manifold $\Si$, such that each fibre
$\pi^{-1}(\si)$ is a real affine straight line in~$\C^3$.

Quite a lot is already known about ruled SL 3-folds in $\C^3$. In 
particular, given a minimal surface $X$ in $\R^3$, Harvey and Lawson 
\cite[\S III.3.C]{HaLa} constructed a ruled SL 3-fold in $\C^3$ 
from the {\it normal bundle} $\nu(X)$ of $X$ in $\R^3$. Borisenko 
\cite[\S 3]{Bori} showed how to generalize this using a harmonic
function $\rho$ on $X$, to define special Lagrangian {\it twisted 
normal bundles} in~$\C^3$.

An important family of ruled special Lagrangian 3-folds are the {\it 
special Lagrangian cones} $N_0$ in $\C^3$. Bryant \cite[Ex.~4]{Brya2} 
showed how to generalize these to the {\it twisted special Lagrangian 
cones} in $\C^3$, using a function $\rho$ on $\Si=N_0\cap{\cal S}^5$
which is an eigenfunction of the Laplacian $\De$ with eigenvalue 2.
Bryant's construction is similar to Borisenko's, but not a generalization.
Bryant also proved other results on ruled SL 3-folds~\cite[\S 7]{Brya2}.

We shall (locally) write each ruled 3-fold $N$ in $\C^3$ in the form
\begin{equation*}
N=\bigl\{r\,\phi(\si)+\psi(\si):r\in\R,\quad \si\in\Si\bigr\},
\end{equation*}
where $\Si$ is a surface, $\phi:\Si\ra{\cal S}^5$ and $\psi:\Si\ra\C^3$
are smooth maps, and ${\cal S}^5$ is the unit sphere in $\C^3$. To each 
ruled 3-fold $N$ we associate an {\it asymptotic cone}
\begin{equation*}
N_0=\bigl\{r\,\phi(\si):r\in\R,\quad \si\in\Si\bigr\}
\end{equation*}
in $\C^3$, to which $N$ is asymptotic at infinity (in a fairly weak sense).

In this paper we shall study the set of ruled SL 3-folds $N$ asymptotic
to a fixed SL cone $N_0$. We begin in \S\ref{ru2} and \S\ref{ru3} by
introducing special Lagrangian geometry, ruled submanifolds and cones,
and \S\ref{ru4} reviews the work of Harvey and Lawson, Borisenko and
Bryant referred to above. The new material begins in \S\ref{ru5},
where we study the equations on $\phi$ and $\psi$ for $N$ to be
special Lagrangian.

It turns out that if $N$ is special Lagrangian then
$N_0$ is, and $\phi$ satisfies a certain nonlinear equation. Taking $N_0$ 
to be special Lagrangian and $\phi$ to be fixed, our first main result, 
Theorem \ref{ru5thm2}, is that provided $N$ is not locally isomorphic to
some $\R^3$ in $\C^3$, it is special Lagrangian if and only if $\psi$
satisfies a {\it linear} equation. Therefore, the set of ruled SL
3-folds $N$ asymptotic to a fixed SL cone $N_0$ in $\C^3$ has the
structure of a vector space.

Our second main result, Theorem \ref{ru6thm1}, is similar to Borisenko's
construction of twisted SL normal bundles, and Bryant's construction of
twisted SL cones. We show that if $N_0$ is a special Lagrangian cone in
$\C^3$ then $\Si$ has the structure of a Riemann surface, and that for 
every holomorphic vector field $w$ on $\Si$ we can construct a ruled SL 
3-fold $N$ asymptotic to $N_0$. Borisenko's result can be regarded as a 
special case of this. Bryant's result is different, and can be combined 
with it to give a larger family of ruled SL 3-folds.

The rest of the paper gives applications of Theorem \ref{ru6thm1}.
In \S\ref{ru6} we show that if $N_0$ is a special Lagrangian cone
on $T^2$ then there is a 2-parameter family of ruled SL 3-folds $N$
asymptotic to $N_0$, and diffeomorphic to $T^2\t\R$. There is also
a variant of this yielding a 1-parameter family of ruled SL 3-folds
diffeomorphic to a nontrivial real line bundle over the Klein bottle.

Section \ref{ru7} gives explicit examples of ruled SL 3-folds.
Using a $\U(1)^2$-invariant $T^2$-cone due to Harvey and Lawson, 
we write down two explicit families of ruled SL 3-folds in $\C^3$
depending on a holomorphic function on $\C$, which include new
kinds of singularities of SL 3-folds. We also use explicit formulae
for a family of SL $T^2$-cones in \cite{Joyc1} to give an explicit
family of ruled SL 3-folds diffeomorphic to~$T^2\t\R$.
\medskip

\noindent{\it Acknowledgements.} I would like to thank Robert Bryant
for helpful conversations, and for bringing reference \cite{Bori} to
my attention.

\section{Special Lagrangian submanifolds in $\C^m$}
\label{ru2}

We begin by defining {\it calibrations} and {\it calibrated 
submanifolds}, following Harvey and Lawson~\cite{HaLa}.

\begin{dfn} Let $(M,g)$ be a Riemannian manifold. An {\it oriented
tangent $k$-plane} $V$ on $M$ is a vector subspace $V$ of
some tangent space $T_xM$ to $M$ with $\dim V=k$, equipped
with an orientation. If $V$ is an oriented tangent $k$-plane
on $M$ then $g\vert_V$ is a Euclidean metric on $V$, so 
combining $g\vert_V$ with the orientation on $V$ gives a 
natural {\it volume form} $\vol_V$ on $V$, which is a 
$k$-form on~$V$.

Now let $\vp$ be a closed $k$-form on $M$. We say that
$\vp$ is a {\it calibration} on $M$ if for every oriented
$k$-plane $V$ on $M$ we have $\vp\vert_V\le \vol_V$. Here
$\vp\vert_V=\al\cdot\vol_V$ for some $\al\in\R$, and 
$\vp\vert_V\le\vol_V$ if $\al\le 1$. Let $N$ be an 
oriented submanifold of $M$ with dimension $k$. Then 
each tangent space $T_xN$ for $x\in N$ is an oriented
tangent $k$-plane. We say that $N$ is a {\it calibrated 
submanifold} if $\vp\vert_{T_xN}=\vol_{T_xN}$ for all~$x\in N$.
\label{ru2def1}
\end{dfn}

It is easy to show that calibrated submanifolds are automatically
{\it minimal submanifolds} \cite[Th.~II.4.2]{HaLa}. Here is the 
definition of special Lagrangian submanifolds in $\C^m$, taken
from~\cite[\S III]{HaLa}.

\begin{dfn} Let $\C^m$ have complex coordinates $(z_1,\dots,z_m)$, 
and define a metric $g$, a real 2-form $\om$ and a complex $m$-form 
$\Om$ on $\C^m$ by
\e
\begin{split}
g=\ms{\d z_1}+\cdots+\ms{\d z_m},\quad
\om&=\frac{i}{2}(\d z_1\w\d\bar z_1+\cdots+\d z_m\w\d\bar z_m),\\
\text{and}\quad\Om&=\d z_1\w\cdots\w\d z_m.
\end{split}
\label{ru2eq1}
\e
Then $\Re\Om$ and $\Im\Om$ are real $m$-forms on $\C^m$. Let
$L$ be an oriented real submanifold of $\C^m$ of real dimension 
$m$, and let $\th\in[0,2\pi)$. We say that $L$ is a {\it special 
Lagrangian submanifold} of $\C^m$, with {\it phase} ${\rm e}^{i\th}$,
if $L$ is calibrated with respect to $\cos\th\Re\Om+\sin\th\Im\Om$, 
in the sense of Definition~\ref{ru2def1}. 

We will often abbreviate `special Lagrangian' by `SL', and 
`$m$-dimensional submanifold' by `$m$-fold', so that we shall talk 
about SL $m$-folds in $\C^m$. Usually we take $\th=0$, so that $L$
has phase 1, and is calibrated with respect to $\Re\Om$. When we
discuss special Lagrangian submanifolds without specifying a phase,
we mean them to have phase~1.
\label{ru2def2}
\end{dfn}

Harvey and Lawson \cite[Cor.~III.1.11]{HaLa} give the following 
alternative characterization of special Lagrangian submanifolds.

\begin{prop} Let\/ $L$ be a real\/ $m$-dimensional submanifold of\/ $\C^m$. 
Then $L$ admits an orientation making it into an SL submanifold of\/ 
$\C^m$ with phase ${\rm e}^{i\th}$ if and only if\/ $\om\vert_L\equiv 0$ 
and\/~$(\sin\th\Re\Om-\cos\th\Im\Om)\vert_L\equiv 0$.
\label{ru2prop1}
\end{prop}

Note that an $m$-dimensional submanifold $L$ in $\C^m$ is called 
{\it Lagrangian} if $\om\vert_L\equiv 0$. Thus special Lagrangian 
submanifolds with phase ${\rm e}^{i\th}$ are Lagrangian submanifolds 
satisfying the extra condition that $(\sin\th\Re\Om-\cos\th\Im\Om)
\vert_L\equiv 0$, which is how they get their name.

\section{Ruled submanifolds of $\C^m$ and cones}
\label{ru3}

We now set up some notation for discussing {\it ruled}\/ submanifolds 
in $\C^m$, which will be used in the rest of the paper.

\begin{dfn} Let $N$ be a real $k$-dimensional submanifold in $\C^m$. A 
{\it ruling} $(\Si,\pi)$ of $N$ is a $(k\!-\!1)$-dimensional manifold 
$\Si$ and a smooth map $\pi:N\ra\Si$, such that for each $\si\in\Si$ 
the fibre $\pi^{-1}(\si)$ is a real affine straight line in $\C^m$. 
A {\it ruled submanifold} is a triple $(N,\Si,\pi)$, where $N$ is
a submanifold of $\C^m$ and $(\Si,\pi)$ a ruling of~$N$.

Usually we will refer to the ruled submanifold as $N$, taking $\Si,\pi$
to be given. An {\it r-orientation} for $(\Si,\pi)$ is a choice of
orientation for the real line $\pi^{-1}(\si)$ for each $\si\in\Si$,
which varies continuously with $\si$. A ruled submanifold $(N,\Si,\pi)$
with an r-orientation is called an {\it r-oriented ruled submanifold}.

Let $(N,\Si,\pi)$ be an r-oriented ruled submanifold of $N$, and 
let ${\cal S}^{2m-1}$ be the unit sphere in $\C^m$. Define a map 
$\phi:\Si\ra{\cal S}^{2m-1}$ such that $\phi(\si)$ is the unique 
unit vector parallel to $\pi^{-1}(\si)$ and in the positive 
direction with respect to the orientation on $\pi^{-1}(\si)$, for 
each $\si\in\Si$. It is easy to see that $\phi$ is a smooth map.

Define a map $\psi:\Si\ra\C^m$ such that $\psi(\si)$ is the unique
vector in $\pi^{-1}(\si)$ orthogonal to $\phi(\si)$, for each 
$\si\in\Si$. Then $\psi$ is smooth, and we have
\e
N=\bigl\{r\,\phi(\si)+\psi(\si):\si\in\Si,\quad r\in\R\bigr\}.
\label{ru3eq1}
\e
\label{ru3def1}
\end{dfn}

Note that some submanifolds $N$ of $\C^m$, such as vector subspaces
$\R^k$ for $k\ge 2$, may admit more than one distinct ruling $(\Si,\pi)$,
$(\Si',\pi')$. In this case we consider $(N,\Si,\pi)$ and $(N,\Si',\pi')$
to be different ruled submanifolds. Though we generally refer to a ruled
submanifold as $N$, dropping $\Si,\pi$, we always have a particular ruling
in mind. Also, whether a ruled submanifold admits an r-orientation is
essentially independent of whether it admits an orientation.

Some explanation of what we mean by `submanifold' is called for here.
Sometimes we will treat submanifolds as {\it embedded}, and regard them
as subsets of $\C^m$. But mostly we allow $N$ to be an {\it immersed}\/
submanifold. That is, $N$ is a real $k$-dimensional manifold together
with an immersion $\iota:N\ra\C^m$, which need not be an embedding.
We generally suppress the immersion~$\iota$. 

However, we do not identify $N$ with its image in $\C^m$. In 
particular, two points $p,q\in N$ with $\iota(p)=\iota(q)$ can have
$\pi(p)\ne\pi(q)$, so that the domain of $\pi$ is $N$, not its image
in $\C^m$. Also, locally we can always choose an r-orientation for a
ruling $(\Si,\pi)$, even though globally an r-orientation may not
exist, as we will see in examples later.

Next we discuss {\it cones} in~$\C^m$.

\begin{dfn} A (singular) submanifold $N$ in $\C^m$ is called a {\it cone}, 
with vertex 0, if whenever $p\in N$ then $r\,p\in N$ for all $r\ge 0$. 
Let $N$ be a cone in $\C^m$. We call $N$ {\it two-sided}\/ if $N=-N$, or 
equivalently if whenever $p\in N$ then $r\,p\in N$ for all $r\in\R$. If 
$N$ is not two-sided, we call it {\it one-sided}.
\label{ru3def2}
\end{dfn}

Cones $N$ in $\C^m$, other than vector subspaces $\R^k$, are 
{\it singular} at their vertex 0. So we will have to deal with 
singular submanifolds. In discussing singularities, there is a 
tension between the embedded and immersed points of view. When 
dealing with embedded submanifolds we regard $N$ as a subset of 
$\C^m$, and a singular point as a point in $\C^m$ where $N$ does
not satisfy the usual submanifold conditions.

In the immersed case, in general we should regard $N$ as a singular
manifold together with an immersion $\iota:N\ra\C^m$. However, there
is a class of singularities of immersed submanifolds which arise when
$N$ is a {\it nonsingular} manifold, and $\iota:N\ra\C^m$ is a smooth 
map which is not an immersion at every point. A singular point is then
the image $q=\iota(p)$ in $\C^m$ of a point $p\in N$ where $\iota$ is 
not an immersion. Note that $\iota^{-1}(q)$ may have positive dimension.

We can use this point of view to interpret two-sided cones as 
examples of ruled submanifolds.

\begin{dfn} Let $N_0$ be a $k$-dimensional, two-sided cone in $\C^m$, 
with an isolated singularity at 0, and regarded for the moment
simply as a subset of $\C^m$. Then $N_0\cap{\cal S}^{2m-1}$ is a
nonsingular $(k\!-\!1)$-dimensional submanifold of ${\cal S}^{2m-1}$,
closed under the action of $-1:{\cal S}^{2m-1}\ra{\cal S}^{2m-1}$.
Define $\Si=(N_0\cap{\cal S}^{2m-1})/\{\pm 1\}$. Then $\Si$ is a
nonsingular $(k\!-\!1)$-manifold.

Define a subset $\tilde N_0$ of $\Si\t\C^m$ by $\tilde N_0=\bigl\{(\{\pm p\},rp):
p\in N_0\cap{\cal S}^{2m-1}$, $r\in\R\bigr\}$. Then $\tilde N_0$ is a
nonsingular $k$-manifold. Define $\iota:\tilde N_0\ra\C^m$ by $\iota:
(\{\pm p\},rp)\mapsto rp$ and $\pi:\tilde N_0\ra\Si$ by $\iota:(\{\pm p\},rp)
\mapsto\{\pm p\}$. Then $\iota(\tilde N_0)=N_0$. Also, $\iota$ is an immersion
except on $\iota^{-1}(0)\cong\Si$. Therefore, we may consider $\tilde N_0$ to
be a singular immersed submanifold of $\C^m$, with immersion~$\iota$. 

Clearly, $(\tilde N_0,\Si,\pi)$ is a {\it ruled submanifold}, with image 
$N_0$ in $\C^m$. Thus, every two-sided cone $N_0$ in $\C^m$ may
be regarded as a ruled submanifold. We will generally suppress the 
notation $\tilde N_0$. We shall call $N_0$ an {\it r-oriented two-sided 
cone} if we are given an r-orientation for~$(\tilde N_0,\Si,\pi)$.
\label{ru3def3}
\end{dfn}

Let $N_0$ be an r-oriented two-sided cone. Then as in Definition 
\ref{ru3def1} we can define maps $\phi:\Si\ra{\cal S}^{2m-1}$ and 
$\psi:\Si\ra\C^m$ such that $N_0$ is written in the form \eq{ru3eq1}. 
From the definition we see that $\psi\equiv 0$. Conversely, any ruled
submanifold of the form \eq{ru3eq1} with $\psi\equiv 0$ is obviously
an r-oriented two-sided cone.

\begin{dfn} Let $(N,\Si,\pi)$ be a ruled $k$-dimensional submanifold 
in $\C^m$. Define
\e
N_0=\bigl\{{\bf v}\in\C^m:\text{$\bf v$ is parallel to $\pi^{-1}(\si)$
for some $\si\in\Si$}\bigr\}.
\label{ru3eq2}
\e
Then $N_0$ is usually a $k$-dimensional two-sided cone in $\C^m$.
We call it the {\it asymptotic cone} of~$N$.
\label{ru3def4}
\end{dfn}

The statement that $N_0$ is a $k$-dimensional cone calls for caution
here. For instance, it may be that $k>1$, but that all the fibres of
$\pi:N\ra\Si$ are parallel, so that $N_0$ is just a straight line in
$\C^m$, of dimension~1. 

What we can say is that if $N$ is an immersed submanifold with 
immersion $\iota:N\ra\C^m$, then $N_0$ is the image of a smooth map 
$\iota_0:N\ra\C^m$, but $\iota_0$ may not be an immersion at any point. 
In particular, if $N$ is r-oriented then we may identify $N$ with 
$\Si\t\R$ as a manifold and define $\iota:N\ra\C^m$ by $\iota:(\si,r)\ra 
r\,\phi(\si)+\psi(\si)$, as in \eq{ru3eq1}, and $\iota_0:N\ra\C^m$ is 
given by~$\iota_0:(\si,r)\ra r\,\phi(\si)$.

To explain in what sense a ruled submanifold $N$ is asymptotic to 
its asymptotic cone $N_0$, we make the following definition.

\begin{dfn} Let $N_0$ be a closed cone in $\C^m$, nonsingular except at 
0, and let $N$ be a closed, nonsingular submanifold in $\C^m$. We say 
that $N$ is {\it asymptotic to $N_0$ with order} $O(r^\al)$ for some 
$\al<1$ if there exists a compact subset $K$ in $N$, a constant $R>0$ 
and a diffeomorphism $\Phi:N_0\sm \,\ov{\!B}_R(0)\ra N\sm K$ such that
\e
\begin{split}
&\bmd{\Phi({\bf x})-{\bf x}}=O(r^\al),\quad
\bmd{\nabla\Phi-I}=O(r^{\al-1})\quad\text{and}\\
&\bmd{\nabla^k\Phi}=O(r^{\al-k}) \quad\text{for $k=2,3,\dots$,
as $r\ra\iy$.}
\end{split}
\label{ru3eq3}
\e
Here $\,\ov{\!B}_R(0)$ is the closed ball of radius $R$ in $\C^m$,
$r$ is the radius function on $\C^m$, and $I$ is the identity map 
on $\C^m$. If $N$ is asymptotic to some cone $N_0$ with order 
$O(r^\al)$ for some $\al<1$ then we say $N$ is {\it asymptotically 
conical} with order~$O(r^\al)$.
\label{ru3def5}
\end{dfn}

For a ruled submanifold $N$ with asymptotic cone $N_0$ written in the form
\begin{equation*}
N=\bigl\{r\,\phi(\si)+\psi(\si):\si\in\Si,\quad r\in\R\bigr\}
\quad\text{and}\quad
N_0=\bigl\{r\,\phi(\si):\si\in\Si,\quad r\in\R\bigr\},
\end{equation*}
we may define $\Phi$ by $\Phi:r\,\phi(\si)\mapsto r\,\phi(\si)+\psi(\si)$
for $\si\in\Si$ and $\md{r}>R$. If $\Si$ is compact, so that $\psi$ is
bounded, it is not difficult to show that $N$ is asymptotic to $N_0$ with
order~$O(1)$. 

This is quite a weak form of convergence, as it says only that $N$ 
stays within a fixed distance of $N_0$ near infinity in $\C^m$, rather 
than converging to $N_0$. However, some ruled submanifolds converge 
more strongly than this, and in \S\ref{ru61} we will describe a class 
of ruled SL 3-folds which are asymptotic to their asymptotic cones
with order~$O(r^{-1})$.

Here is an elementary result on ruled special Lagrangian $m$-folds 
in $\C^m$. The proof is easy, and we omit it.

\begin{prop} Let\/ $N$ be a ruled special Lagrangian $m$-fold in $\C^m$,
and let\/ $N_0$ be the asymptotic cone of\/ $N$, as in Definition
\ref{ru3def4}. Then $N_0$ is a special Lagrangian cone in $\C^m$,
provided it is $m$-dimensional.
\label{ru3prop}
\end{prop}

The reason for requiring $N_0$ to be $m$-dimensional here is 
that there do exist ruled SL $m$-folds $N$ whose tangent cones
have $\dim N_0<m$, and so cannot be regarded as special Lagrangian
except in a rather singular sense. For instance, if $N=\Si\t\R$ in 
$\C^{m-1}\t\C$, where $\Si$ is an SL $(m\!-\!1)$-fold in $\C^{m-1}$, 
then all the lines in the obvious ruling of $N$ are parallel, and so
$N_0$ is just $\R$ in~$\C^m$.

Motivated by this proposition, the general point of view we will take 
is to fix a special Lagrangian cone $N_0$ in $\C^3$, and study the 
ruled special Lagrangian 3-folds $N$ in $\C^3$ asymptotic to $N_0$. 
In the notation of Definition \ref{ru2def1}, $N_0$ determines the 
map $\phi:\Si\ra{\cal S}^5$, and we shall look for maps
$\psi:\Si\ra\C^3$ such that $N$ defined by \eq{ru2eq1} is special
Lagrangian.

\section{Review of previous work}
\label{ru4}

Before beginning our new material, we briefly review three previous
papers that have contributed to the theory of ruled special Lagrangian
3-folds. These are Harvey and Lawson \cite{HaLa}, Borisenko \cite{Bori}
and Bryant~\cite{Brya2}.

\subsection{Harvey and Lawson's SL normal bundles in $\C^m$}
\label{ru41}

Harvey and Lawson \cite[\S III.3.C]{HaLa} gave the following
construction of SL $m$-folds in $\C^m$. Let $\R^m$ have coordinates
$(x_1,\ldots,x_m)$ and Euclidean metric $\d x_1^2+\cdots+\d x_m^2$,
let $X$ be a submanifold of $\R^m$, and let $\nu(X)$ be the normal
bundle of $X$. That is, for each $x\in X$ the fibre $\nu_x$ of $\nu(X)$
is the orthogonal complement $T_xX^\perp$ of $T_xX$ in~$\R^m$.

Write $\nu(X)$ as a subset of $\R^m\op\R^m$ by
\begin{equation*}
\nu=\bigl\{(x,y):x\in X,\quad y\in\nu_x=T_xX^\perp\subset\R^n\bigr\}.
\end{equation*}
Then it is a classical fact that $\nu(X)$ is Lagrangian with respect 
to the symplectic structure $\om=\d x_1\w\d y_1+\cdots+\d x_m\w\d y_m$,
where $(x_1,\ldots,x_m,y_1,\ldots,y_m)$ are the obvious coordinates 
on $\R^m\op\R^m$. Identify $\R^m\op\R^m$ with $\C^m$ using the complex 
coordinates $(z_1,\ldots,z_m)$, where $z_j=x_j+iy_j$. Then $\om$ is 
the usual K\"ahler form on $\C^m$. Thus, for any submanifold $X$ of $\R^m$, 
the normal bundle $\nu(X)$ is a Lagrangian submanifold of~$\C^m$.

Harvey and Lawson were interested in the conditions on $X$ for $\nu(X)$ to 
be special Lagrangian. Their answer is given in the following definition 
\cite[Def.~III.3.15]{HaLa} and theorem~\cite[Th.~III.3.11]{HaLa}.

\begin{dfn} Let $X$ be a $k$-dimensional submanifold of $\R^m$, 
and $A$ the second fundamental form of $X$. Then $A_x$ lies in 
$\nu_x\ot S^2T_x^*X$ for each $x\in X$. We call $X$ {\it austere}
if for all $x\in X$ and $y\in\nu_x$, the invariants of odd order of 
the quadratic form $A_x\cdot y$ in $S^2T_x^*X$ vanish. Equivalently, 
for all $x\in X$ and $y\in\nu_x$, the collection of eigenvalues 
$\la_1,\dots,\la_k$ of $A_x\cdot y$, with multiplicity, should 
be invariant under multiplication by~$-1$.
\label{ru4def1}
\end{dfn}

Austere submanifolds of $\R^n$ are studied by Bryant~\cite{Brya1}.

\begin{thm} Let\/ $X$ be a $k$-dimensional submanifold of\/ $\R^m$.
Then the normal bundle $\nu(X)$ is special Lagrangian in $\C^m$,
with phase $i^{m-k}$, if and only if\/ $X$ is austere.
\label{ru4thm1}
\end{thm}

As the sign of the phase does not matter, we can take $\nu(X)$ to have
phase 1 if $m-k$ is even, and phase $i$ if $m-k$ is odd.

Now the condition for $X$ to be {\it minimal}\/ in $\R^n$ is that
$\Tr(A_x\cdot y)=0$ for all $x\in X$ and $y\in\nu_x$, as $\Tr(A)$ 
is the mean curvature of $X$. But $\Tr(A_x\cdot y)$ is an odd-order
invariant of $A_x\cdot y$, and so vanishes when $X$ is austere.
Thus all austere submanifolds are minimal. Furthermore, when $X$ is 
2-dimensional, $\Tr(A_x\cdot y)$ is the {\it only} odd-order invariant
of $A_x\cdot y$, and so $X$ is austere if and only if it is minimal.
So we prove:

\begin{cor} Let $X$ be a minimal surface in $\R^m$. Then the normal 
bundle $\nu(X)$ is special Lagrangian in $\C^m$, with phase~$i^{m-2}$.
\label{ru4cor1}
\end{cor}

We shall be concerned only with the case $m=3$. When $k=\dim X$ is 0,1 
or 3, the construction above yields only affine special Lagrangian 
3-planes $\R^3$ in $\C^3$. But when $k=2$, by considering minimal 
surfaces in $\R^3$ it yields nontrivial examples of special Lagrangian 
3-folds $\nu(X)$ in $\C^3$, with phase $i$. Note that in this case
$\nu(X)$ is automatically a {\it ruled}\/ special Lagrangian 3-fold, 
with ruling $(X,\pi)$, where $\pi:\nu(X)\ra X$ is the natural projection,
whose fibres are real straight lines in~$\C^3$.

\subsection{Borisenko's twisted SL normal bundles in $\C^3$}
\label{ru42}

For each minimal surface $X$ in $\R^3$, Harvey and Lawson's
construction yields a normal bundle $\nu(X)$ in $\C^3$, which
is a ruled SL 3-fold. In \cite[\S 3]{Bori}, Borisenko generalized 
this construction to generate a family of ruled SL 3-folds in $\C^3$,
which we will call {\it twisted special Lagrangian normal bundles},
depending on a minimal surface $X$ in $\R^3$ and a harmonic function
$\rho:X\ra\R$. Here is Borisenko's result~\cite[Th.~1]{Bori}.

\begin{thm} Let\/ $X$ be an oriented, regular minimal surface in
$\R^3$ and $\rho:X\ra\R$ a harmonic function. Let\/ $(s,t)$ be
local coordinates on $X$ compatible with the orientation, and 
write the immersion $X\ra\R^3$ as $(s,t)\mapsto{\bf x}(s,t)$.
Define vector-valued functions ${\bf n},{\bf p}:X\ra\R^3$ by
\e
{\bf n}=\frac{{\bf x}_s\t{\bf x}_t}{\md{{\bf x}_s\t{\bf x}_t}}
\quad\text{and}\quad
{\bf p}=\frac{\rho_s{\bf x}_t-\rho_t{\bf x}_s}{
\md{{\bf x}_s\t{\bf x}_t}}\t{\bf n},
\label{ru4eq1}
\e
where ${\bf x}_s=\frac{\pd{\bf x}}{\pd s}$, and so on. Then 
${\bf n},{\bf p}$ are well-defined and independent of the 
choice of coordinates $(s,t)$. Define
\e
N=\bigl\{{\bf x}+i\bigl({\bf p}({\bf x})+r\,{\bf n}({\bf x})\bigr):
{\bf x}\in X,\quad r\in\R\bigr\}.
\label{ru4eq2}
\e
Then $N$ is a ruled special Lagrangian $3$-fold in $\C^3=\R^3\op i\,\R^3$ 
with phase~$i$.
\label{ru4thm2}
\end{thm}

Here `$\t$' is the usual cross product on $\R^3$. It is easy to
see that $\bf n$ is the positive unit normal to $X$, defined using
the orientations on $X$ and $\R^3$, and $\bf p$ is the gradient
vector of $\rho$, that is ${\bf p}^a=g^{ab}(\d\rho)_b$ in index
notation, where $g$ is the metric on $X$ induced from the Euclidean
metric on~$\R^3$.

If we reverse the orientation on $X$ then $\bf n$ changes sign and
$\bf p$ is fixed. Reversing the sign of $t$ in \eq{ru4eq2}, we see
that $N$ is unchanged by reversing the orientation of $X$. In fact
the construction works for non-orientable $X$ as well.

When $\rho$ is constant, Borisenko's construction reduces to
Harvey and Lawson's special Lagrangian normal bundle $\nu(X)$. 
If $N$ is one of Borisenko's twisted normal bundles, then the
{\it asymptotic cone} $N_0$ of $N$ is a subset of the special 
Lagrangian 3-plane $i\,\R^3$ in~$\C^3$.

\subsection{Bryant's results on ruled SL 3-folds}
\label{ru43}

In \cite{Brya2}, Bryant proves a number of results on
ruled special Lagrangian 3-folds. Here is one, given in 
\cite[Ex.~4]{Brya2} and based on work in \cite[\S 4]{Brya1},
which is similar to Borisenko's result.

\begin{thm} Let\/ $N_0$ be an r-oriented, two-sided special
Lagrangian cone in $\C^3$, with ruling $(\Si,\pi)$. Then $\Si$ is
an oriented Riemannian $2$-manifold with an isometric immersion
$\phi:\Si\ra{\cal S}^5$. Suppose $\Si$ is simply-connected. 
Let\/ $\rho:\Si\ra\R$ be any solution of the second-order, 
linear elliptic equation~$*\d(*\d\rho)+2\rho=0$. 

Define a $\C^3$-valued\/ $1$-form $\be$ on $\Si$ by $\be=\phi*\d\rho-
\rho *\d\phi$. Then $\be$ is closed, so there exists ${\bf b}:\Si\ra\C^3$ 
with\/ $\d{\bf b}=\be$. Define $N=\bigl\{r\,\phi(\si)+{\bf b}(\si):
\si\in\Si$, $r\in\R\bigr\}$. Then $N$ is a ruled special Lagrangian 
$3$-fold in $\C^3$, asymptotic to~$N_0$.
\label{ru4thm3}
\end{thm}

When $\rho\equiv 0$ we have $\be=0$ and $\bf b$ is constant, so $N$ 
is a translation of $N_0$ in $\C^3$. Thus these examples generalize
special Lagrangian cones in $\C^3$, and Bryant calls them {\it twisted
special Lagrangian cones}. They do not agree with Borisenko's twisted
special Lagrangian normal bundles.

Next we summarize the results of \cite[\S 3.7]{Brya2}, which considers
the family of all ruled SL 3-folds in $\C^3$. Let $\La$ be the set
of oriented lines in $\C^3$. Then $\La$ is a 10-dimensional manifold,
which fibres over ${\cal S}^5$ with fibre $\R^5$. Bryant shows that
$\La$ carries a {\it real analytic, Levi-flat almost CR structure} 
$(E,J)$, which is a subbundle $E$ of $T\La$ with fibre $\R^8$ (in 
this case), and an almost complex structure $J$ on the fibres of 
$E$, satisfying some conditions.

A ruled 3-fold $N$ in $\C^3$ may be regarded as the total space 
of a 2-dimensional family of (oriented) real lines, and thus as 
a 2-manifold $\Si$ in $\La$. Bryant shows that $N$ is special
Lagrangian if and only if $\Si$ is $E$-{\it holomorphic}, that 
is, the tangent spaces of $\Si$ lie in $E$ and are closed under~$J$.

Bryant also uses the theory of exterior differential systems to 
measure how `big' the various families of ruled SL 3-folds are. 
In the sense of exterior differential systems, Harvey and Lawson's 
family of special Lagrangian normal bundles in $\C^3$, and the family 
of special Lagrangian cones in $\C^3$, both depend on two functions 
of one variable. 

Similarly, Borisenko's family of twisted special Lagrangian normal 
bundles in $\C^3$, and Bryant's family of twisted special Lagrangian 
normal bundles in $\C^3$, both depend on four functions of one variable.
But the family of all ruled special Lagrangian 3-folds in $\C^3$ depends
on {\it six} functions of one variable. So generic ruled SL 3-folds do 
not come from Theorems~\ref{ru4thm1}--\ref{ru4thm3}.

\section{General results on ruled SL 3-folds}
\label{ru5}

In this section we shall parametrize an r-oriented ruled special 
Lagrangian 3-fold $N$ in $\C^3$ using maps $\phi:\Si\ra{\cal S}^5$ 
and $\psi:\Si\ra\C^3$ as in Definition \ref{ru3def1}, and determine 
the conditions on $\phi$ and $\psi$ for $N$ to be special Lagrangian. 

The following notation will be used throughout this section. Let 
$\Si$ be a 2-dimensional, connected, real analytic manifold. Let 
$\phi:\Si\ra{\cal S}^5$ be a real analytic immersion, where 
${\cal S}^5$ is the unit sphere in $\C^3$. Let $\psi:\Si\ra\C^3$ 
be a real analytic map. Define
\e
N=\bigl\{r\,\phi(\si)+\psi(\si):\si\in\Si,\quad r\in\R\bigr\}.
\label{ru5eq1}
\e
Then $N$ is an r-oriented ruled 3-fold in $\C^3$. We suppose $N$ 
is special Lagrangian.

Now $N$ is the image of the real analytic map $\Phi:\Si\t\R\ra\C^3$ 
given by
\e
\Phi(\si,r)=r\,\phi(\si)+\psi(\si).
\label{ru5eq2}
\e
As $\phi$ is an immersion, $\Phi$ is an immersion almost everywhere
in $\Si\t\R$. The images of points where $\Phi$ is not an immersion
are generally singular points of $N$. Regarding $N$ as an immersed
copy of $\Si\t\R$ with (possibly singular) immersion $\Phi$, we may
define $\pi:N\ra\Si$ by $\pi:(\si,r)\ra\si$, and then $(\Si,\pi)$ is
a ruling of~$N$.

The asymptotic cone $N_0$ of $N$ is 
\e
N_0=\bigl\{r\,\phi(\si):\si\in\Si,\quad r\in\R\bigr\},
\label{ru5eq3}
\e
which is the image of $\Phi_0:\Si\t\R\ra\C^3$ given by
\e
\Phi_0(\si,r)=r\,\phi(\si).
\label{ru5eq4}
\e
As we have assumed $\phi$ is an immersion, $\Phi_0$ is an
immersion except when $r=0$, so that $N_0$ is nonsingular
as an immersed submanifold except at 0. Note that $N_0$ is
special Lagrangian, by Proposition~\ref{ru3prop}.

As $\phi:\Si\ra{\cal S}^5$ is an immersion, the pull-back
$\phi^*(g)$ of the round metric on ${\cal S}^5$ makes $\Si$
into a {\it Riemannian $2$-manifold}. Also, as $N_0$ is special
Lagrangian it is oriented, since $\Re\Om$ is a nonvanishing 3-form 
on $N_0$. We use this to define a natural {\it orientation} on 
$\Si$, such that the non-vanishing 2-form $\phi^*(\phi\cdot\Re\Om)$
on $\Si$ is positive. Equivalently, local coordinates $(s,t)$ on 
$\Si$ are oriented if
\e
{\ts\Re\Om\bigl(\phi,\frac{\pd\phi}{\pd s},\frac{\pd\phi}{\pd t}\bigr)>0.} 
\label{ru5eq5}
\e

Thus $\Si$ is an oriented Riemannian 2-manifold. Therefore it has 
a natural {\it complex structure} $J$. Near any point $\si\in\Si$ 
we can choose a holomorphic coordinate $z=s+it$. The corresponding
real coordinates $(s,t)$ on $\Si$ have the property that
$J\bigl(\frac{\pd}{\pd s}\bigr)\equiv\frac{\pd}{\pd t}$. We shall
call local coordinates $(s,t)$ on $\Si$ with this property 
{\it oriented conformal coordinates}.

Section \ref{ru51} analyzes the conditions on $\phi$ and $\psi$ for
$N$ to be special Lagrangian, and \S\ref{ru52} studies ruled SL 3-folds 
using an `evolution equation' approach. These ideas are combined in 
\S\ref{ru53} to prove our first main result, Theorem \ref{ru5thm2}, 
which will be applied in the rest of the paper.

\subsection{The special Lagrangian equations on $\phi$ and $\psi$}
\label{ru51}

We shall find the conditions on $\phi,\psi$ for $N$ to be 
special Lagrangian at each nonsingular point. Suppose $\Phi$ 
is an immersion at $(\si,r)$ in $\Si\t\R$, and let $p=\Phi(\si,r)$. 
Choose oriented conformal coordinates $(s,t)$ on $\Si$ near $\si$. 
Then $T_pN=\an{{\bf v}_1,{\bf v}_2,{\bf v}_3}_{\R}$, where
\e
\begin{gathered}
{\bf v}_1=\frac{\pd\Phi}{\pd r}(\si,r)=\phi(\si),\qquad
{\bf v}_2=\frac{\pd\Phi}{\pd s}(\si,r)=r\frac{\pd\phi}{\pd s}(\si)
+\frac{\pd\psi}{\pd s}(\si)\\
\text{and}\qquad
{\bf v}_3=\frac{\pd\Phi}{\pd t}(\si,r)=r\frac{\pd\phi}{\pd t}(\si)+
\frac{\pd\psi}{\pd t}(\si).
\end{gathered}
\label{ru5eq6}
\e

We need $T_pN$ to be a special Lagrangian 3-plane $\R^3$ in $\C^3$,
with phase 1. By Proposition \ref{ru2prop1}, the condition for this 
is that $\om\vert_{T_pN}=\Im\Om\vert_{T_pN}=0$, which is equivalent to 
\e
\om({\bf v}_1,{\bf v}_2)=\om({\bf v}_1,{\bf v}_3)=
\om({\bf v}_2,{\bf v}_3)=\Im\Om({\bf v}_1,{\bf v}_2,{\bf v}_3)=0.
\label{ru5eq7}
\e
Substituting in for ${\bf v}_1,{\bf v}_2,{\bf v}_3$ using \eq{ru5eq6}
gives equations upon $\phi$ and $\psi$ and their derivatives, which
are linear or quadratic polynomials in $r$. As the equations should
hold for all $r\in\R$, the coefficient of each power of $r$ should 
vanish. 

So we find that \eq{ru5eq7} holding for all $r$ is equivalent to
the equations
\ea
{\ts\om\bigl(\phi,\frac{\pd\phi}{\pd s}\Bigr)=
\om\bigl(\phi,\frac{\pd\phi}{\pd t}\bigr)=
\om\bigl(\frac{\pd\phi}{\pd s},\frac{\pd\phi}{\pd t}\bigr)=
\Im\Om\bigl(\phi,\frac{\pd\phi}{\pd s},\frac{\pd\phi}{\pd t}\bigr)}&=0,
\label{ru5eq8}\\
\begin{split}
{\ts \om\bigl(\phi,\frac{\pd\psi}{\pd s}\Bigr)=
\om\bigl(\phi,\frac{\pd\psi}{\pd t}\bigr)=
\om\bigl(\frac{\pd\phi}{\pd s},\frac{\pd\psi}{\pd t}\bigr)+
\om\bigl(\frac{\pd\psi}{\pd s},\frac{\pd\phi}{\pd t}\bigr)}&=0,\\
{\ts\Im\Om\bigl(\phi,\frac{\pd\phi}{\pd s},\frac{\pd\psi}{\pd t}\bigr)+
\Im\Om\bigl(\phi,\frac{\pd\psi}{\pd s},\frac{\pd\phi}{\pd t}\bigr)}&=0,
\end{split}
\label{ru5eq9}\\
\text{and}\qquad
{\ts\om\bigl(\frac{\pd\psi}{\pd s},\frac{\pd\psi}{\pd t}\bigr)=
\Im\Om\bigl(\phi,\frac{\pd\psi}{\pd s},\frac{\pd\psi}{\pd t}\bigr)}&=0.
\label{ru5eq10}
\ea
Here we have arranged the equations so that \eq{ru5eq8} does 
not involve $\psi$ at all, \eq{ru5eq9} is linear in $\psi$, and 
\eq{ru5eq10} is quadratic in $\psi$. Note that \eq{ru5eq8} implies 
that $\an{\phi,\frac{\pd\phi}{\pd s},\frac{\pd\phi}{\pd t}}_{\R}$ 
is an SL 3-plane in $\R^3$, which is the condition for $N_0$ to be 
special Lagrangian. 

In the next proposition we shall show that \eq{ru5eq9} and \eq{ru5eq10} 
are equivalent to one of two linear equations on $\psi$. We will need the 
following notation, adapted from \cite[\S 5]{Joyc3} with an extra factor 
of 2. Define an anti-bilinear cross product $\t:\C^3\t\C^3\ra\C^3$ by
\e
({\bf u}\t{\bf v})^b=2{\bf u}^{a_1}{\bf v}^{a_2}(\Re\Om)_{a_1a_2a_3}
g^{a_3b},
\label{ru5eq11}
\e
using the index notation for (real) tensors on $\C^3$. Calculation
using \eq{ru2eq1} shows that in coordinates `$\t$' is given by
\e
(u_1,u_2,u_3)\t(v_1,v_2,v_3)=(\bar u_2\bar v_3-\bar u_3\bar v_2,
\bar u_3\bar v_1-\bar u_1\bar v_3,\bar u_1\bar v_2-\bar u_2\bar v_1).
\label{ru5eq12}
\e
It is equivariant under the $\SU(3)$-action, as $\Re\Om$ and $g$ are
$\SU(3)$-invariant.

\begin{prop} In the situation above, $N$ is special Lagrangian
if and only~if
\e
\om\Bigl(\phi,\frac{\pd\phi}{\pd s}\Bigr)\equiv 0
\quad\text{and}\quad
\frac{\pd\phi}{\pd t}\equiv\phi\t\frac{\pd\phi}{\pd s},
\label{ru5eq13}
\e
and\/ $\psi$ satisfies either
\begin{itemize}
\item[{\rm(i)}] $\om\bigl(\phi,\frac{\pd\psi}{\pd s}\bigr)\equiv 0$ and\/
$\frac{\pd\psi}{\pd t}\equiv\phi\t\frac{\pd\psi}{\pd s}+f\phi$, where 
`$\t$' is defined in \eq{ru5eq12} and\/ $f:\Si\ra\R$ is some real 
function; or
\item[{\rm(ii)}] $\frac{\pd\psi}{\pd s}(\si)$ and\/ 
$\frac{\pd\psi}{\pd t}(\si)$ lie in $\ban{\phi(\si),\frac{\pd\phi}{\pd s}
(\si),\frac{\pd\phi}{\pd t}(\si)}_{\R}$ for all\/~$\si\in\Si$.
\end{itemize}
\label{ru5prop1}
\end{prop}

\begin{proof} Above we showed that $N$ is special Lagrangian if and only 
if \eq{ru5eq8}--\eq{ru5eq10} hold. We will show that \eq{ru5eq13} is 
equivalent to \eq{ru5eq8}, and \eq{ru5eq9}--\eq{ru5eq10} are equivalent 
to (i) or (ii). Fix $\si\in\Si$, and let $C=\bmd{\frac{\pd}{\pd s}(\si)}$. 
Then $C>0$. As $(s,t)$ are oriented conformal coordinates and $\md{\phi}
\equiv 1$, one can show that if \eq{ru5eq8} holds then $\phi(\si),
C^{-1}\frac{\pd\phi}{\pd s}(\si),C^{-1}\frac{\pd\phi}{\pd t}(\si)$ are 
an oriented orthonormal basis of a special Lagrangian 3-plane in~$\C^3$. 

Let $(w_1,w_2,w_3)$ be the unique complex coordinate system on $\C^3$ 
in which we have $\phi(\si)=(1,0,0)$, $C^{-1}\frac{\pd\phi}{\pd s}(\si)=
(0,1,0)$ and $C^{-1}\frac{\pd\phi}{\pd t}(\si)=(0,0,1)$. Then 
$(w_1,w_2,w_3)$ are related to the usual coordinates $(z_1,z_2,z_3)$
on $\C^3$ by an $\SU(3)$ transformation. Thus, as $g,\om$ and $\Om$
are $\SU(3)$-invariant, by \eq{ru2eq1} we have
\e
\begin{split}
g=\ms{\d w_1}+\cdots+\ms{\d w_3},\quad
\om&=\frac{i}{2}(\d w_1\w\d\bar w_1+\cdots+\d w_3\w\d\bar w_3),\\
\text{and}\quad\Om&=\d w_1\w\d w_2\w\d z_3.
\end{split}
\label{ru5eq14}
\e

In the coordinates $(w_1,w_2,w_3)$, write
\ea
\phi(\si)=(1,0,0),\qquad
{\ts\frac{\pd\phi}{\pd s}}(\si)&=(0,C,0),\qquad
{\ts\frac{\pd\phi}{\pd t}}(\si)=(0,0,C),
\label{ru5eq15}\\
\begin{split}
{\ts\frac{\pd\psi}{\pd s}}(\si)&=(a_1+ib_1,a_2+ib_2,a_3+ib_3),\\
\text{and}\qquad
{\ts\frac{\pd\psi}{\pd t}}(\si)&=(c_1+id_1,c_2+id_2,c_3+id_3),
\end{split}
\label{ru5eq16}
\ea
for $a_j,b_j,c_j,d_j\in\R$. Using \eq{ru5eq14}--\eq{ru5eq16}, 
equations \eq{ru5eq9} and \eq{ru5eq10} become
\ea
b_1=d_1=d_2-b_3=d_3+b_2&=0,
\label{ru5eq17}\\
a_1d_1+a_2d_2+a_3d_3-b_1c_1-b_2c_2-b_3c_3&=0,
\label{ru5eq18}\\
\text{and}\qquad
a_2d_3+b_2c_3-a_3d_2-b_3c_2&=0,
\label{ru5eq19}
\ea
cancelling factors of~$C$.

Setting $b_1=d_1=0$ and substituting $d_2=b_3$ and $d_3=-b_2$
by \eq{ru5eq17}, equations \eq{ru5eq18} and \eq{ru5eq19} give
\begin{equation*}
-b_2(a_3+c_2)+b_3(a_2-c_3)=0 \quad\text{and}\quad
-b_2(a_2-c_3)-b_3(a_3+c_2)=0,
\end{equation*}
which is equivalent to
\e
\begin{pmatrix} -b_2 & b_3 \\ -b_3 & -b_2 \end{pmatrix}
\begin{pmatrix} a_3+c_2 \\ a_2-c_3 \end{pmatrix}=
\begin{pmatrix} 0 \\ 0 \end{pmatrix}.
\label{ru5eq20}
\e

The $2\t 2$ matrix appearing here has determinant $b_2^2+b_3^2$. If 
this is nonzero then the matrix is invertible, so the column matrix
must be zero. If the determinant is zero then $b_2=b_3=0$. So 
\eq{ru5eq20} holds if and only if either
\begin{itemize}
\item[(a)] $a_3+c_2=a_2-c_3=0$, or
\item[(b)] $b_2=b_3=0$.
\end{itemize}
These two possibilities correspond to parts (i) and (ii) of the proposition.

As the transformation from the $(z_1,z_2,z_3)$ coordinates to 
the $(w_1,w_2,w_3)$ coordinates lies in $\SU(3)$, and `$\t$' is 
$\SU(3)$-equivariant, the formula for `$\t$' in the $(w_1,w_2,w_3)$ 
coordinates is the same as \eq{ru5eq12} in the $(z_1,z_2,z_3)$ 
coordinates. Thus, combining \eq{ru5eq12} and \eq{ru5eq15} we see
that $\frac{\pd\phi}{\pd t}(\si)=\phi(\si)\t\frac{\pd\phi}{\pd s}(\si)$.
Since this holds for all $\si\in\Si$, equation \eq{ru5eq8} implies
\eq{ru5eq13}. Conversely, it is easy to show that \eq{ru5eq13} 
implies~\eq{ru5eq8}.

Suppose (a) holds, so that $c_2=-a_3$ and $c_3=a_2$. Then by 
\eq{ru5eq17} we have
\e
{\ts\frac{\pd\psi}{\pd s}}(\si)=(a_1,a_2\!+\!ib_2,a_3\!+\!ib_3)
\quad\text{and}\quad
{\ts\frac{\pd\psi}{\pd t}}(\si)=(c_1,-a_3\!+\!ib_3,a_2\!-\!ib_2).
\label{ru5eq21}
\e
Clearly $\om\bigl(\phi(\si),\frac{\pd\psi}{\pd s}(\si)\bigr)=0$ as
$b_1=0$. Equations \eq{ru5eq12}, \eq{ru5eq15} and \eq{ru5eq21} give
\begin{equation*}
{\ts\frac{\pd\psi}{\pd t}}(\si)=\phi(\si)\t{\ts\frac{\pd\psi}{\pd s}}(\si)
+c_1\phi(\si).
\end{equation*}
Thus part (i) of the proposition holds at $\si$, with $f(\si)=c_1$. 
Conversely, if (i) holds at $\si$ then \eq{ru5eq21} holds, so 
\eq{ru5eq17}--\eq{ru5eq19} hold, and thus \eq{ru5eq9} and \eq{ru5eq10} 
hold at~$\si$.

Now suppose (b) holds, so that $b_2=b_3=0$. Then \eq{ru5eq17} gives
\begin{equation*}
{\ts\frac{\pd\psi}{\pd s}}(\si)=(a_1,a_2,a_3)
\quad\text{and}\quad
{\ts\frac{\pd\psi}{\pd t}}(\si)=(c_1,c_2,c_3).
\end{equation*}
That is, $\frac{\pd\psi}{\pd s}(\si)$ and $\frac{\pd\psi}{\pd t}(\si)$
are real vectors in the $(w_1,w_2,w_3)$ coordinates. It follows that
part (ii) of the proposition holds at $\si$. Conversely, if part (ii) 
holds at $\si$ then \eq{ru5eq17}--\eq{ru5eq19} hold, and so \eq{ru5eq9} 
and \eq{ru5eq10} hold at~$\si$.

We have shown that \eq{ru5eq13} is equivalent to \eq{ru5eq8}, and
\eq{ru5eq9} and \eq{ru5eq10} hold if and only if (i) or (ii) holds at 
every point $\si$ in $\Si$. These are not exclusive options; both (i) 
and (ii) may hold at some points. But we need one of (i) and (ii) to 
hold at {\it every} point in $\Si$, rather than (i) at some points and (ii) at others. 

Let $\si$ be a generic point in $\Si$. Then if (i) or (ii) holds at $\si$, 
it also holds in a small neighbourhood of $\si$ in $\Si$. But $\Si$ is 
connected, $\phi,\psi$ are real analytic, and (i), (ii) are closed 
conditions, so if one of (i) or (ii) holds in an open set in $\Si$ 
then it holds in all of $\Si$. This completes the proof.
\end{proof}

Observe that parts (i) and (ii) are both {\it linear} restrictions on $\psi$.
So $\psi$ satisfies one of two linear equations. The proposition suggests
that there are really two different kinds of ruled special Lagrangian 
3-folds. However, we will show later that the only ruled SL 3-folds 
admitting rulings satisfying (ii) but not (i) are SL 3-planes $\R^3$ 
in $\C^3$, so that all interesting ruled SL 3-folds satisfy part~(i).

Note also that $N$ is unchanged by transformations of the form
\e
\phi\mapsto\phi, \qquad \psi\mapsto \psi+\al\phi,
\label{ru5eq22}
\e
where $\al:\Si\ra\R$ is a real analytic function. This can be regarded
as a kind of {\it gauge transformation} of $(\phi,\psi)$. We can fix
$\al$ and $\psi$ uniquely by requiring that $g(\phi,\psi)\equiv 0$,
as we did in Definition \ref{ru3def1}, but we will not always do this.
Alternatively, we can use $\al$ to ensure that the function $f$ in
part (i) of the proposition is zero.

\subsection{Evolution equations for ruled SL 3-folds}
\label{ru52}

We will now study ruled special Lagrangian 3-folds using 
the `evolution equation' approach developed by the author in 
\cite{Joyc1,Joyc2,Joyc3}. This depends on the following 
result, proved in~\cite[Th.~3.3]{Joyc1}.

\begin{thm} Let\/ $P$ be a compact, orientable, real analytic 
$(m-1)$-manifold, $\chi$ a real analytic, nonvanishing section 
of\/ $\La^{m-1}TP$, and\/ $\Psi:P\ra\C^m$ a real analytic
immersion such that\/ $\Psi^*(\om)\equiv 0$ on $P$. Then there 
exists $\ep>0$ and a unique real analytic family 
$\bigl\{\Psi_t:t\in(-\ep,\ep)\bigr\}$ of real analytic maps 
$\Psi_t:P\ra\C^m$ with $\Psi_0=\Psi$, satisfying the equation
\e
\Bigl(\frac{\d\Psi_t}{\d t}\Bigr)^b=(\Psi_t)_*(\chi)^{a_1\ldots a_{m-1}}
(\Re\Om)_{a_1\ldots a_{m-1}a_m}g^{a_mb},
\label{ru5eq23}
\e
using the index notation for (real) tensors on $\C^m$. Define 
$\Phi:P\t(-\ep,\ep)\ra\C^m$ by $\Phi(p,t)=\Psi_t(p)$. Then 
$N=\Im\Phi$ is a nonsingular immersed special Lagrangian 
submanifold of\/~$\C^m$.
\label{ru5thm1}
\end{thm}

Here the assumption that $P$ is compact is often unnecessary. The 
theorem constructs SL $m$-folds in $\C^m$ by evolving arbitrary
real analytic $(m\!-\!1)$-submanifolds $P$ in $\C^m$ with
$\om\vert_P\equiv 0$. The trouble with this is that as the set of
such submanifolds is infinite-dimensional, the theorem is really
an infinite-dimensional evolution problem, and so difficult to
solve explicitly.

To get round this, in \cite{Joyc1,Joyc2,Joyc3} the author found
various special classes $\cal C$ of real analytic $(m-1)$-submanifolds
$P$ in $\C^m$ with $\om\vert_P\equiv 0$, such that the evolution
\eq{ru5eq23} stayed within the class $\cal C$. Usually $\cal C$ was
only finite-dimensional, so that \eq{ru5eq23} reduced to an o.d.e.\
in finitely many variables, which might even be solved explicitly.

We will use the same idea to study ruled SL 3-folds in $\C^3$ in our 
next proposition. We shall see that if $\cal C$ is the class of 
ruled 2-manifolds $P$ in $\C^3$ with $\om\vert_P\equiv 0$, then the 
evolution \eq{ru5eq23} stays within $\cal C$, and the resulting SL 
3-folds $N$ are ruled. In what follows, we say that a function
defined on a compact interval $S$ in $\R$ is {\it real analytic} if
it extends to a real analytic function on an open neighbourhood of 
$S$ in~$\R$.

\begin{prop} Let\/ $S$ be the circle $\R/\Z$, or a compact interval in 
$\R$, and let\/ $s$ be a coordinate on $S$, taking values in $\R/\Z$ 
or $\R$. Suppose $\phi_0:S\ra{\cal S}^5$ and\/ $\psi_0:S\ra\C^3$ are 
real analytic maps satisfying
\e
\om\Bigl(\phi_0,\frac{\pd\phi_0}{\pd s}\Bigr)\equiv
\om\Bigl(\phi_0,\frac{\pd\psi_0}{\pd s}\Bigr)\equiv 0
\quad\text{in $S$.}
\label{ru5eq24}
\e
Then there exists $\ep>0$ and unique real analytic maps 
$\phi:S\t(-\ep,\ep)\ra{\cal S}^5$ and\/ $\psi:S\t(-\ep,\ep)\ra\C^3$ 
with\/ $\phi(s,0)=\phi_0(s)$ and\/ $\psi(s,0)=\psi_0(s)$ for all\/
$s\in S$, satisfying
\begin{gather}
\frac{\pd\phi}{\pd t}=\phi\t\frac{\pd\phi}{\pd s},\qquad
\frac{\pd\psi}{\pd t}=\phi\t\frac{\pd\psi}{\pd s}
\label{ru5eq25}\\
\text{and}\quad
\om\Bigl(\phi,\frac{\pd\phi}{\pd s}\Bigr)\equiv
\om\Bigl(\phi,\frac{\pd\psi}{\pd s}\Bigr)\equiv 0
\quad\text{in $S\t(-\ep,\ep)$.}
\label{ru5eq26}
\end{gather}
Define $N=\bigl\{\Phi(r,s,t):r\in\R$, $s\in S$, $t\in(-\ep,\ep)\bigr\}$,
where $\Phi:(r,s,t)\mapsto r\,\phi(s,t)+\psi(s,t)$. Then $N$ is an
r-oriented ruled special Lagrangian $3$-fold in~$\C^3$. 
\label{ru5prop2}
\end{prop}

\begin{proof} We shall apply Theorem \ref{ru5thm1}. Define $P$ to be 
$\R\t S$, with coordinates $(r,s)$, and let $\chi=2\frac{\pd}{\pd r}
\w\frac{\pd}{\pd s}$. Then $\chi$ is a real analytic, nonvanishing 
section of $\La^2TP$, as in the theorem.

Consider a family of maps $\Psi_t:P\ra\C^3$ depending smoothly on 
$t\in(-\ep,\ep)$, of the form $\Psi_t(r,s)=r\,\phi(s,t)+\psi(s,t)$, 
where $\phi,\psi:S\t(-\ep,\ep)\ra\C^3$ are real analytic maps with 
$\phi(s,0)\equiv\phi_0(s)$ and $\psi(s,0)\equiv\psi_0(s)$. We do
not yet assume that $\phi$ maps to ${\cal S}^5$ in $\C^3$, but we 
will prove this later.

Then $\Psi_t^*\bigl(\frac{\pd}{\pd r}\bigr)=\phi$ and $\Psi_t^*\bigl(
\frac{\pd}{\pd s}\bigr)=r\frac{\pd\phi}{\pd s}+\frac{\pd\psi}{\pd s}$, 
so \eq{ru5eq23} becomes
\begin{align*}
r\frac{\pd\phi^b}{\pd t}+\frac{\pd\psi^b}{\pd t}&=
2\phi^{a_1}\Bigl(r\frac{\pd\phi^{a_2}}{\pd s}+
\frac{\pd\psi^{a_2}}{\pd s}\Bigr)(\Re\Om)_{a_1a_2a_3}g^{a_3b}\\
&=\Bigl(r\,\phi\t\frac{\pd\phi}{\pd s}+
\phi\t\frac{\pd\psi}{\pd s}\Bigr)^b,
\end{align*}
by the definition of `$\t$' in \eq{ru5eq11}. Equating linear and
constant terms in $r$ on each side, we see that \eq{ru5eq23} 
is equivalent to~\eq{ru5eq25}.

The only problem in applying Theorem \ref{ru5thm1} in this
situation is that $P$ is not compact. However, as the $\Psi_t$
are defined using maps $\phi,\psi$ on $S$, which is compact, this 
does not matter, and the theorem applies. This can be checked by
examining its proof in~\cite[\S 3]{Joyc1}.

Thus, by Theorem \ref{ru5thm1}, there exists $\ep>0$ and unique
real analytic maps $\phi,\psi:S\t(-\ep,\ep)\ra\C^3$ with 
$\phi(s,0)\equiv\phi_0(s)$ and $\psi(s,0)\equiv\psi_0(s)$,
satisfying \eq{ru5eq25}, and such that $N$ is special Lagrangian.
It remains only to prove \eq{ru5eq26}, and that $\phi$ maps 
into~${\cal S}^5$. 

Equation \eq{ru5eq26} follows from the equation $\Phi^*(\om)\bigl(
\frac{\pd}{\pd r},\frac{\pd}{\pd s}\bigr)=0$, which holds as $N$ is 
Lagrangian. And $\frac{\pd}{\pd t}\bigl(\ms{\phi}\bigr)=
2g\bigl(\phi,\frac{\pd\phi}{\pd t}\bigr)=
2g\bigl(\phi,\phi\t\frac{\pd\phi}{\pd s}\bigr)=0$,
by \eq{ru5eq25} and the definition of `$\t$'. Thus $\md{\phi}$ is
independent of $t$. But $\bmd{\phi(s,0)}=1$ as $\phi(s,0)=\phi_0(s)$
lies in ${\cal S}^5$, so $\md{\phi}\equiv 1$, and $\phi$ 
maps~$S\t(-\ep,\ep)\ra{\cal S}^5$.
\end{proof}

Observe the similarity between equation \eq{ru5eq13} and part (i) of 
Proposition \ref{ru5prop1}, and equations \eq{ru5eq25} and \eq{ru5eq26}
of Proposition \ref{ru5prop2}. The only real difference is that
part (i) of Proposition \ref{ru5prop1} allows $\frac{\pd\psi}{\pd t}
\equiv\phi\t\frac{\pd\psi}{\pd s}+f\phi$ for some $f:\Si\ra\R$, whereas
\eq{ru5eq25} requires $f$ to be zero. We reconcile these by noting, as
above, that we can set $f$ to be zero using a transformation of the
form \eq{ru5eq22} depending on $\al:\Si\ra\R$. The `evolution' 
construction has the effect of fixing~$\al$.
\medskip

\noindent{\bf Remark.} One consequence of the proposition is that if
$P$ is a real analytic ruled 2-manifold in $\C^3$ with $\om\vert_P=0$,
then $P$ extends locally to a unique ruled special Lagrangian 3-fold
$N$ in $\C^3$. This fact is implicit in Bryant's Cartan--K\"ahler 
theory calculations \cite[\S 3.7]{Brya2}, and so is not new.

\subsection{Main results}
\label{ru53}

The next proposition deals with ruled SL 3-folds satisfying part (ii)
but not part (i) of Proposition~\ref{ru5prop1}.

\begin{prop} Every ruled special Lagrangian $3$-fold\/ $N$ in $\C^3$
locally admits an r-oriented ruling $(\Si,\pi)$ satisfying part\/ 
{\rm(i)} of Proposition \ref{ru5prop1}. An r-oriented ruled special 
Lagrangian $3$-fold\/ $(N,\Si,\pi)$ in $\C^3$ satisfying part\/ {\rm(ii)} 
but not part\/ {\rm(i)} of Proposition \ref{ru5prop1} is locally 
isomorphic to an affine special Lagrangian $3$-plane $\R^3$ in~$\C^3$.
\label{ru5prop3}
\end{prop}

\begin{proof} Suppose $N$ is a ruled SL 3-fold in $\C^3$, with r-oriented
ruling $(\Si,\pi)$ satisfying part (ii) but not part (i) of Proposition 
\ref{ru5prop1}. As $N$ is real analytic wherever is is nonsingular by 
\cite[Th.~III.2.7]{HaLa}, we may take $(\Si,\pi)$ to be real analytic,
at least locally.

Choose any embedded real analytic curve $\ga:[0,1]\ra\Si$. Define 
$S=[0,1]$ and $\phi_0:S\ra{\cal S}^5$, $\psi_0:S\ra\C^3$ by
$\phi_0(s)=\phi(\ga(s))$ and $\psi_0(s)=\psi(\ga(s))$. Apply
Proposition \ref{ru5prop2}. This gives real analytic maps
$\phi':S\t(-\ep,\ep)$ and $\psi':S\t(-\ep,\ep)\ra\C^3$,
and constructs a ruled SL 3-fold $N'$ in $\C^3$ from them,
with a ruling $(\Si',\pi')$ that satisfies part (i) of
Proposition~\ref{ru5prop1}.

Now $N$ and $N'$ intersect in the ruled surface $\pi^{-1}(\Im\ga)$
in $\C^3$. Therefore, by \cite[Th.~III.5.5]{HaLa}, $N$ and $N'$
coincide locally. Thus, locally $N$ admits a ruling $(\Si',\pi')$
satisfying part (i) of Proposition \ref{ru5prop1}. Since by assumption
$(\Si,\pi)$ does not satisfy (i), $(\Si,\pi)$ and $(\Si',\pi')$ must
be different.

It is not difficult to show that different, nearby curves $\ga$ in $\Si$
will yield different rulings $(\Si',\pi')$ of $N$. Thus $N$ admits not
just two, but infinitely many different rulings. But Bryant 
\cite[Th.~6]{Brya2} shows that any SL 3-fold with more than two distinct
rulings is planar, that is, locally isomorphic to $\R^3$ in~$\C^3$.
\end{proof}

Here is our main result, which follows from 
Propositions~\ref{ru5prop1}--\ref{ru5prop3}.

\begin{thm} Let\/ $(N,\Si,\pi)$ be a non-planar, r-oriented, ruled 
SL\/ $3$-fold in $\C^3$. Then there exist real analytic maps 
$\phi:\Si\ra{\cal S}^5$ and\/ $\psi:\Si\ra\C^3$ such that
\e
N=\bigl\{r\,\phi(\si)+\psi(\si):\si\in\Si,\quad r\in\R\bigr\}.
\label{ru5eq27}
\e
If\/ $(s,t)$ are oriented conformal coordinates on $U\subset\Si$, then 
$\phi,\psi$ satisfy
\begin{gather}
\om\Bigl(\phi,\frac{\pd\phi}{\pd s}\Bigr)\equiv 
\om\Bigl(\phi,\frac{\pd\psi}{\pd s}\Bigr)\equiv 0,
\label{ru5eq28}\\
\frac{\pd\phi}{\pd t}=\phi\t\frac{\pd\phi}{\pd s}
\quad\text{and}\quad
\frac{\pd\psi}{\pd t}\equiv\phi\t\frac{\pd\psi}{\pd s}+f\phi
\label{ru5eq29}
\end{gather}
for some real analytic function $f:U\ra\R$. Conversely, if\/
$\phi,\psi$ satisfy these equations then $N$ is special Lagrangian
wherever it is nonsingular.
\label{ru5thm2}
\end{thm}

Again, note that \eq{ru5eq28} and \eq{ru5eq29} are {\it linear} in 
$\psi$, regarding $\phi$ as fixed and $f$ as linear in $\psi$. This
means that the family of ruled special Lagrangian 3-folds $N$ with
a fixed asymptotic cone $N_0$ has the structure of a vector space.

\section{Holomorphic vector fields and ruled SL 3-folds}
\label{ru6}

In Theorems \ref{ru4thm2} and \ref{ru4thm3} we described constructions 
of ruled SL 3-folds in $\C^3$ by Borisenko and Bryant. Borisenko's result 
involved a harmonic function $\rho$ on a minimal surface $X$ in $\R^3$, 
and Bryant's involved an SL cone $N_0$ in $\C^3$ and a function $\rho$ 
on $\Si=N_0\cap{\cal S}^5$ satisfying~$*\d(*\d\rho)+2\rho=0$.

We shall now present a construction of ruled special Lagrangian
3-folds in $\C^3$ which is similar to both Borisenko's and
Bryant's constructions, but not the same as either. The data
we use is an SL cone $N_0$ in $\C^3$ and a holomorphic vector
field $w$ on the Riemann surface $\Si=N_0\cap{\cal S}^5$. Here 
is our result.

\begin{thm} Let\/ $N_0$ be an r-oriented, two-sided special 
Lagrangian cone in $\C^3$. Then as in \S\ref{ru5} we can write
\e
N_0=\bigl\{r\,\phi(\si):\si\in\Si,\quad r\in\R\bigr\},
\label{ru6eq1}
\e
where $\Si$ is a Riemann surface and $\phi:\Si\ra{\cal S}^5$ a
real analytic map, such that if\/ $(s,t)$ are oriented conformal
coordinates on an open set\/ $U\subset\Si$ then $\phi$ satisfies
\e
\om\Bigl(\phi,\frac{\pd\phi}{\pd s}\Bigr)\equiv 0
\quad\text{and}\quad
\frac{\pd\phi}{\pd t}=\phi\t\frac{\pd\phi}{\pd s}.
\label{ru6eq2}
\e

Let\/ $w$ be a holomorphic vector field on $\Si$ and define
$\psi:\Si\ra\C^3$ by $\psi={\cal L}_w\phi$, where ${\cal L}_w$
is the Lie derivative. Define
\e
N=\bigl\{r\,\phi(\si)+\psi(\si):\si\in\Si,\quad r\in\R\bigr\}.
\label{ru6eq3}
\e
Then $N$ is an r-oriented ruled special Lagrangian $3$-fold in~$\C^3$. 
\label{ru6thm1}
\end{thm}

\begin{proof} The first part of the theorem follows from the 
material of \S\ref{ru5}. If $w\equiv 0$ then the second part
is trivial. So suppose $w\not\equiv 0$, so that $w$ has only 
isolated zeros. Let $\si\in\Si$ be a point where $w$ is nonzero. 
Then it is easy to show that as $w$ is a holomorphic vector field 
on $\Si$, there exist oriented conformal coordinates $(s,t)$ on an 
open neighbourhood $U$ of $\si$ in $\Si$ such that $w=\frac{\pd}{\pd s}$ 
in $U$. Therefore $\psi=\frac{\pd\phi}{\pd s}$ in~$U$.

Taking $\frac{\pd}{\pd s}$ of the first equation of \eq{ru6eq2}, we find
\e
\om\Bigl(\frac{\pd\phi}{\pd s},\frac{\pd\phi}{\pd s}\Bigr)+
\om\Bigl(\phi,\frac{\pd^2\phi}{(\pd s)^2}\Bigr)\equiv
\om\Bigl(\phi,\frac{\pd\psi}{\pd s}\Bigr)\equiv 0,
\label{ru6eq4}
\e
as $\om$ is antisymmetric and $\psi=\frac{\pd\phi}{\pd s}$. Similarly,
taking $\frac{\pd}{\pd s}$ of the second equation of \eq{ru6eq2} gives
\begin{equation*}
\frac{\pd^2\phi}{\pd s\pd t}=
\frac{\pd\phi}{\pd s}\t\frac{\pd\phi}{\pd s}+
\phi\t\frac{\pd^2\phi}{(\pd s)^2}=
\phi\t\frac{\pd^2\phi}{(\pd s)^2},
\end{equation*}
as `$\t$' is antisymmetric. Using $\frac{\pd^2}{\pd s\pd t}=
\frac{\pd^2}{\pd t\pd s}$ and $\psi=\frac{\pd\phi}{\pd s}$,
this becomes
\e
\frac{\pd\psi}{\pd t}=\phi\t\frac{\pd\psi}{\pd s}.
\label{ru6eq5}
\e

Comparing \eq{ru6eq2}, \eq{ru6eq4} and \eq{ru6eq5} with equations
\eq{ru5eq28} and \eq{ru5eq29}, we see by Theorem \ref{ru5thm2} that
the subset $\pi^{-1}(U)$ of $N$ is a ruled SL 3-fold. Now every point 
$\si\in\Si$ has such a neighbourhood $U$, except for the isolated 
zeros of $w$. Thus $N$ is an r-oriented ruled SL 3-fold, except 
possibly along a discrete set of lines. But to be special Lagrangian 
is a closed condition on the nonsingular part of $N$, so $N$ is 
special Lagrangian.
\end{proof}

Here is a sketch of how the theorem relates to Borisenko's 
construction in Theorem \ref{ru4thm2}. Let $X$ be an oriented 
minimal surface in $\R^3$, with unit normal ${\bf n}:X\ra{\cal S}^2$.
Then $\bf n$ is the {\it Gauss map} of $X$. Suppose for simplicity
that $\bf n$ is a diffeomorphism between $X$ and an open set $\Si$
of~${\cal S}^2$. 

Note that $\bf n$ is a conformal map, though it is some sense 
orientation-reversing. Thus the harmonic function $\rho:X\ra\R$ 
pushes forward to a harmonic function $\rho':\Si\ra\R$ under 
$\bf n$, since $X$ and $\Si$ are conformal. Define $\al=\d\rho'
-iJ\d\rho'$. Then $\al$ is a holomorphic (1,0)-form on~$\Si$.

Now using the second fundamental form of $X$ one can construct a
natural holomorphic section $\be$ of $T^{(1,0)}\Si\ot_{\C} T^{(1,0)}\Si$.
The contraction $\al\cdot\be$ is a holomorphic section of $T^{(1,0)}\Si$.
Let $w$ be the real part of $\al\cdot\be$. Then $w$ is a holomorphic
vector field on $\Si$. The map $\phi:\Si\ra{\cal S}^5$ is just 
multiplication by $i$, and $\psi={\cal L}_w\phi$ in Theorem \ref{ru6thm1} 
corresponds to Borisenko's term $i{\bf p}:\Si\ra\C^3$ in~\eq{ru4eq2}.

Thus $\psi$ in Theorem \ref{ru6thm1} corresponds to the `twisting' term
$\bf p$ in Theorem \ref{ru4thm2}. However, our result works for all SL 
cones in $\C^3$, whereas Borisenko's construction works only when the 
asymptotic cone is $i\,\R^3$. To make up for this, Borisenko also 
includes the extra term $\bf x$ in \eq{ru4eq2}, which does not appear 
in our construction.

Bryant's construction in Theorem \ref{ru4thm3} does not coincide
with our construction above, even though both begin with a special
Lagrangian cone. In fact, because of the linearity in $\psi$ noted
after Theorem \ref{ru5thm2} we can combine the two constructions to get:

\begin{thm} Let\/ $N_0,\Si$ and\/ $\phi$ be as above, and suppose $\Si$ 
is simply-connected. Let\/ $w$ be a holomorphic vector field on $\Si$ 
and define $\psi:\Si\ra\C^3$ by $\psi={\cal L}_w\phi$, where ${\cal L}_w$
is the Lie derivative. Let\/ $\rho:\Si\ra\R$ be any solution of the 
second-order, linear elliptic equation $*\d(*\d\rho)+2\rho=0$. Define a 
$\C^3$-valued\/ $1$-form $\be$ on $\Si$ by $\be=\phi*\d\rho-\rho *\d\phi$. 
Then $\be$ is closed, so there exists ${\bf b}:\Si\ra\C^3$ with\/ 
$\d{\bf b}=\be$. Define
\e
N=\bigl\{r\,\phi(\si)+\psi(\si)+{\bf b}(\si):\si\in\Si,\quad r\in\R\bigr\}.
\label{ru6eq6}
\e
Then $N$ is an r-oriented ruled special Lagrangian $3$-fold in~$\C^3$.
\label{ru6thm2}
\end{thm}

To see that the terms $\psi$ and $\bf b$ here represent genuinely 
different perturbations of $N_0$, observe that in oriented conformal
coordinates $(s,t)$ we have
\begin{alignat*}{2}
{\ts\frac{\pd\psi}{\pd s}}&=
{\ts u\frac{\pd^2\phi}{(\pd s)^2}
+\frac{\pd u}{\pd s}\frac{\pd\phi}{\pd s}
+v\frac{\pd^2\phi}{\pd s\pd t}
+\frac{\pd v}{\pd s}\frac{\pd\phi}{\pd t},}&\quad
{\ts\frac{\pd\psi}{\pd t}}&=
{\ts u\frac{\pd^2\phi}{\pd s\pd t}
+\frac{\pd u}{\pd t}\frac{\pd\phi}{\pd s}
+v\frac{\pd^2\phi}{(\pd t)^2}
+\frac{\pd v}{\pd t}\frac{\pd\phi}{\pd t}},\\
{\ts\frac{\pd{\bf b}}{\pd s}}&={\ts-\frac{\pd\rho}{\pd t}\phi+\rho
\frac{\pd\phi}{\pd t}}\qquad\qquad \text{and} &\quad
{\ts\frac{\pd{\bf b}}{\pd t}}&={\ts\frac{\pd\rho}{\pd s}\phi-\rho
\frac{\pd\phi}{\pd s},}
\end{alignat*}
where~$w=u\frac{\pd}{\pd s}+v\frac{\pd}{\pd t}$. 

Generically five of $\phi,\frac{\pd\phi}{\pd s},\frac{\pd\phi}{\pd t},
\frac{\pd^2\phi}{(\pd s)^2},\frac{\pd^2\phi}{\pd s\pd t}$ and 
$\frac{\pd^2\phi}{(\pd t)^2}$ will be linearly independent, and so 
the $\psi$ and $\bf b$ terms cannot agree, since if they did then
only four would be linearly dependent. In fact, in the generic case 
the author expects that locally Theorem \ref{ru6thm2} gives {\it all}\/ 
the ruled special Lagrangian 3-folds asymptotic to $N_0$. But this will 
not in general be true at points where $\phi$ and its first and second 
derivatives are too linearly dependent.

\subsection{Ruled SL 3-folds over compact Riemann surfaces}
\label{ru61}

We shall now apply Theorem \ref{ru6thm1} in the case when $\Si$ is 
a {\it compact}\/ Riemann surface, without boundary. If $\Si$ is a 
compact, connected Riemann surface of genus $g$ then the vector
space of holomorphic vector fields on $\Si$ has dimension 6 when
$g=0$, dimension 2 when $g=1$, and dimension 0 when $g\ge 2$. So
to get nontrivial holomorphic vector fields we should take $\Si$ 
to be ${\cal S}^2$ or~$T^2$.

However, it follows from well known facts in minimal surface theory
that any special Lagrangian cone on ${\cal S}^2$ in $\C^3$ must be
an SL 3-plane $\R^3$. The author first learnt this from Robert
Bryant, and a proof can be found in Haskins \cite[Th.~B]{hask2}.
When applied to $N_0=\R^3$, Theorem \ref{ru6thm2} just gives back a 
different ruling of the same~$\R^3$.

Therefore the only interesting case is $\Si\cong T^2$. Then we can prove:

\begin{thm} Let\/ $N_0$ be an r-oriented two-sided special Lagrangian 
cone on $T^2$. That is, $N_0$ may be defined as in \eq{ru6eq1}, where
$\Si\cong T^2$ is a Riemann surface and $\phi:\Si\ra{\cal S}^5$ a real 
analytic immersion satisfying \eq{ru6eq2} in oriented conformal 
coordinates. Then there exists a $2$-dimensional family of distinct,
r-oriented, ruled special Lagrangian $3$-folds $N$ with asymptotic 
cone $N_0$, which are asymptotic to $N_0$ with order $O(r^{-1})$ in 
the sense of Definition~\ref{ru3def5}.
\label{ru6thm3}
\end{thm}

\begin{proof} Any Riemann surface $\Si\cong T^2$ may be written
as $\R^2/\La$, where $\La\cong\Z^2$ is a lattice in $\R^2$, and the 
coordinates $(s,t)$ on $\R^2$ are oriented conformal coordinates.
Write $\Si$ in this way. Then the holomorphic vector fields $w$ 
on $\Si$ are of the form $u\frac{\pd}{\pd s}+v\frac{\pd}{\pd t}$
for~$u,v\in\R$. 

Lift $\phi$ to a $\La$-invariant map $\phi:\R^2\ra{\cal S}^5$. 
For each $u,v\in\R$, define
\e
{\ts N_{u,v}=\bigl\{r\,\phi(s,t)+u\frac{\pd\phi}{\pd s}(s,t)
+v\frac{\pd\phi}{\pd t}(s,t):r,s,t\in\R\bigr\}.}
\label{ru6eq7}
\e
Then $N_{u,v}$ is a ruled special Lagrangian 3-fold in $\C^3$
wherever it is nonsingular, by Theorem \ref{ru6thm1}, and it
clearly has asymptotic cone~$N_0$. 

It remains to prove that $N_{u,v}$ is asymptotic to $N_0$ with 
order $O(r^{-1})$, in the sense of Definition \ref{ru3def5}. Let 
$R>0$, and define $\Phi:N_0\sm\,\ov{\!B}_R(0)\ra N_{u,v}$ by
\e
\begin{split}
{\ts\Phi:r\,\phi(s,t)\longmapsto
r\,\phi\bigl(s-\frac{u}{r},t-\frac{v}{r}\bigr)}
&{\ts+u\frac{\pd\phi}{\pd s}\bigl(s-\frac{u}{r},t-\frac{v}{r}\bigr)}\\
&{\ts+v\frac{\pd\phi}{\pd t}\bigl(s-\frac{u}{r},t-\frac{v}{r}\bigr)}
\end{split}
\label{ru6eq8}
\e
for $\md{r}>R$ and $s,t\in\R$. Then $\Phi$ is well-defined, and
using the expansions
\begin{gather*}
{\ts\frac{\pd\phi}{\pd s}\bigl(s-\frac{u}{r},t-\frac{v}{r}\bigr)=
\frac{\pd\phi}{\pd s}(s,t)+O(r^{-1}),\quad
\frac{\pd\phi}{\pd t}\bigl(s-\frac{u}{r},t-\frac{v}{r}\bigr)=
\frac{\pd\phi}{\pd t}(s,t)+O(r^{-1})}\\
\text{and}\quad
{\ts\phi\bigl(s-\frac{u}{r},t-\frac{v}{r}\bigr)=
\phi(s,t)-\frac{u}{r}\frac{\pd\phi}{\pd s}(s,t)
-\frac{v}{r}\frac{\pd\phi}{\pd t}(s,t)+O(r^{-2})}
\end{gather*}
for large $r$, one can show that
\begin{equation*}
\Phi\bigl(r\,\phi(s,t)\bigr)=r\,\phi(s,t)+O(r^{-1})
\end{equation*}
for large $r$, which is the first equation of \eq{ru3eq3}, with 
$\al=-1$. The other equations of \eq{ru3eq3} may be proved in the same 
way. Therefore $N_{u,v}$ is asymptotic to $N_0$ with order~$O(r^{-1})$. 
\end{proof}

In \S\ref{ru3} we saw that ruled submanifolds are asymptotic to 
their asymptotic cones with order $O(1)$. But the ruled SL 3-folds 
in the theorem are asymptotic with order $O(r^{-1})$, which is 
stronger. The reason for this is that each line $\pi^{-1}(\si)$
in $N_0$ is translated by $\psi(\si)$ to make $N$, and $\psi(\si)$
is tangent to $N_0$ along $\pi^{-1}(\si)$. Thus the distance between
$N$ and $N_0$ is roughly proportional to the curvature of $N_0$, which
is~$O(r^{-1})$.

One can show by the same method that any ruled SL 3-fold constructed 
in Theorem \ref{ru6thm1} is asymptotic to $N_0$ with order $O(r^{-1})$ 
over compact subsets of $\Si$. The author believes that all ruled SL 
3-folds asymptotic with order $O(r^{-1})$ to their asymptotic cones 
arise from the construction of Theorem~\ref{ru6thm1}.

In Theorem \ref{ru6thm3} we took $N_0$ to be fibred by a $T^2$ family 
of lines. Thus $N_0\cap{\cal S}^5$ is actually two copies of $T^2$, as 
each line intersects ${\cal S}^5$ in two points, and the two $T^2$ are 
swapped by the action of $-1$ on ${\cal S}^5$. So $N_0$ is two opposite 
$T^2$-cones meeting at their common vertex 0, and for generic $u,v$ we 
expect $N_{u,v}$ to be a nonsingular immersed 3-submanifold diffeomorphic 
to $T^2\t\R$, with two $T^2$ ends at infinity.

However, there is another kind of special Lagrangian $T^2$-cone, in 
which $N_0$ is one $T^2$-cone rather than two, invariant under $\pm 1$. 
Suppose this is the case. Let $\tilde\Si=N_0\cap{\cal S}^5$, so that 
$\tilde\Si$ is a Riemann surface isomorphic to $T^2$. Then the action of 
$-1$ on ${\cal S}^5$ restricts to a free, orientation-reversing involution 
on $\tilde\Si$, and $\Si=\tilde\Si/\{\pm1\}$ is the {\it Klein bottle}.

Now $(N_0,\Si,\pi)$ is a ruled 3-fold as in Definition \ref{ru3def3}. 
However, $N_0$ is {\it not}\/ r-orientable, so that we cannot define
$\phi:\Si\ra{\cal S}^5$ because of sign problems. Instead, we define
$\phi:\tilde\Si\ra{\cal S}^5$ to be the identity map on the double cover
$\tilde\Si$ of $\Si$. This gives a corresponding immersed, r-oriented SL 
$T^2$-cone $\tilde N_0$, isomorphic as a (singular) immersed manifold to
$\tilde\Si\t\R$ with immersion $\iota:(\tilde\si,r)\mapsto r\,\phi(\tilde\si)$,
which is the double cover of~$N_0$.

We may then apply Theorem \ref{ru6thm3} to $\tilde N_0$ to get a 
2-parameter family of ruled SL 3-folds $\tilde N_{u,v}$ asymptotic
to $\tilde N_0$. It turns out that $\tilde N_{u,v}$ is the double cover 
of a ruled SL 3-fold $N_{u,v}$ asymptotic to $N_0$, fibred over 
the Klein bottle, if and only if the holomorphic vector field
$w=u\frac{\pd}{\pd s}+v\frac{\pd}{\pd t}$ changes sign under
$-1:\tilde\Si\ra\tilde\Si$. The vector space of such $w$ is 1-dimensional.
Thus we prove:

\begin{thm} Let\/ $N_0$ be a two-sided special Lagrangian cone on 
the Klein bottle. Then there exists a $1$-dimensional family of 
distinct, non r-orientable, ruled special Lagrangian $3$-folds $N$ 
with asymptotic cone $N_0$, which are asymptotic to $N_0$ with order 
$O(r^{-1})$ in the sense of Definition~\ref{ru3def5}.
\label{ru6thm4}
\end{thm}

When the 3-folds $N$ in the theorem are nonsingular, they are
immersed 3-submanifolds diffeomorphic to the total space of
a nontrivial real line bundle over the Klein bottle. They have
one end at infinity, which is asymptotically a $T^2$-cone.

Finally we note that we can generate interesting ruled SL 3-folds
$N$ from SL cones $N_0$ over Riemann surfaces $\Si$ with any 
genus $g\ge 1$, by applying Theorem \ref{ru6thm1} to nontrivial 
{\it meromorphic} vector fields $w$ on $\Si$ with poles at 
$\si_1,\ldots,\si_k$ in $\Si$. Then $N$ is asymptotic to the
union of $N_0$ and $k$ SL 3-planes $\Pi_1,\ldots,\Pi_k$, which 
are tangent to $N_0$ along the lines $\pi^{-1}(\si_1),\ldots,
\pi^{-1}(\si_k)$ respectively.

\section{Explicit examples of ruled SL 3-folds}
\label{ru7}

We now apply Theorem \ref{ru6thm1} to give some explicit examples
of ruled special Lagrangian 3-folds. We start by considering the 
$\U(1)^2$-invariant SL $T^2$-cone found by Harvey and Lawson 
\cite[\S III.3.A]{HaLa}. Define $\phi:\R^2\ra{\cal S}^5$ by
\e
\phi:(s,t)\longmapsto{\ts\frac{1}{\sqrt{3}}}\bigl({\rm e}^{is},
{\rm e}^{-\frac{i}{2}s-\frac{i\sqrt{3}}{2}t},
{\rm e}^{-\frac{i}{2}s+\frac{i\sqrt{3}}{2}t}\bigr),
\label{ru7eq1}
\e
and write $N_0=\bigl\{r\,\phi(s,t):r,s,t\in\R\bigr\}$. Then
$N_0$ is a special Lagrangian $T^2$-cone, and $(s,t)$ are
oriented conformal coordinates on $\Si=\R^2$, considered
as a Riemann surface.

Note that $\phi$ is invariant under the lattice $\La\cong\Z^2$ in 
$\R^2$ generated by $(2\pi,2\pi/\sqrt{3}\,)$ and $(0,4\pi/\sqrt{3}\,)$,
but we shall not pass to the quotient $\R^2/\La$. Let $u(s,t)+iv(s,t)$ 
be a holomorphic function of $s+it$. Then $w=u(s,t)\frac{\pd}{\pd s}+
v(s,t)\frac{\pd}{\pd t}$ is a holomorphic vector field on $\Si$.
Applying Theorem \ref{ru6thm1} gives: 

\begin{thm} Let\/ $u,v:\R^2\ra\R$ be functions such that\/
$u(s,t)+iv(s,t)$ is a holomorphic function of\/ $s+it$. Define
\e
\begin{split}
N=\Bigl\{{\ts\frac{1}{\sqrt{3}}}\Bigl(
&{\rm e}^{is}\bigl(r+iu(s,t)\bigr),
{\rm e}^{-\frac{i}{2}s-\frac{i\sqrt{3}}{2}t}
{\ts\bigl(r-\frac{i}{2}u(s,t)-\frac{i\sqrt{3}}{2}v(s,t)\bigr)},\\
&{\rm e}^{-\frac{i}{2}s+\frac{i\sqrt{3}}{2}t}
{\ts\bigl(r-\frac{i}{2}u(s,t)+\frac{i\sqrt{3}}{2}v(s,t)\bigr)}
\Bigr):r,s,t\in\R\Bigr\}.
\end{split}
\label{ru7eq2}
\e
Then $N$ is a ruled special Lagrangian $3$-fold in~$\C^3$.
\label{ru7thm1}
\end{thm}

The good thing about this theorem is that it defines a large family
of ruled SL 3-folds very explicitly. Therefore we can use it as a
laboratory for studying the {\it singularities} of ruled SL 3-folds. 
For instance, if we put $u(s,t)+iv(s,t)=(s+it)^k$ for $k=2,3,\dots$,
then we generate a series of ruled SL 3-folds $N_k$ with an isolated 
singularity at 0 in~$\C^3$. 

We can also exchange the dependent and independent variables in
Theorem \ref{ru7thm1}, and regard $s+it$ as a holomorphic function 
of $u+iv$. This yields:

\begin{thm} Let\/ $s,t:\R^2\ra\R$ be functions such that\/
$s(u,v)+it(u,v)$ is a holomorphic function of\/ $u+iv$. Define
\e
\begin{split}
N=\Bigl\{{\ts\frac{1}{\sqrt{3}}}\Bigl(
&{\rm e}^{is(u,v)}\bigl(r+iu\bigr),
{\rm e}^{-\frac{i}{2}s(u,v)-\frac{i\sqrt{3}}{2}t(u,v)}
{\ts\bigl(r-\frac{i}{2}u-\frac{i\sqrt{3}}{2}v\bigr)},\\
&{\rm e}^{-\frac{i}{2}s(u,v)+\frac{i\sqrt{3}}{2}t(u,v)}
{\ts\bigl(r-\frac{i}{2}u+\frac{i\sqrt{3}}{2}v\bigr)}
\Bigr):r,u,v\in\R\Bigr\}.
\end{split}
\label{ru7eq3}
\e
Then $N$ is a ruled special Lagrangian $3$-fold in~$\C^3$.
\label{ru7thm2}
\end{thm}

This is less useful for modelling singularities, as $N$ is always
nonsingular near 0 in $\C^3$. But one can instead use it for modelling
`branched' asymptotic behaviour of ruled SL 3-folds. Putting 
$s(u,v)+it(u,v)=(u+iv)^k$ for $k=2,3,\dots$, we find that near the
line $\bigl\{(x,x,x):x\in\R\bigr\}$, $N$ is a non-singular ruled SL 
3-fold that is asymptotic to a $k$-fold branched cover of the 
Harvey--Lawson $T^2$-cone $N_0$ defined above.

Next we shall generalize a family of $\U(1)$-invariant special Lagrangian 
$T^2$-cones in $\C^3$ given in \cite[\S 8]{Joyc1} to a family of ruled 
SL 3-folds, as illustrations of Theorems \ref{ru6thm3} and \ref{ru6thm4}. 
They are written in terms of the {\it Jacobi elliptic functions}, which 
we now briefly introduce. The following material can be found in 
Chandrasekharan~\cite[Ch.~VII]{Chan}. 

For each $k\in[0,1]$, the Jacobi elliptic functions $\sn(t,k)$, $\cn(t,k)$
and $\dn(t,k)$ with modulus $k$ are the unique solutions to the o.d.e.s
\begin{align*}
\bigl({\textstyle\frac{\d}{\d t}}\sn(t,k)\bigr)^2&=\bigl(1-\sn^2(t,k)\bigr)
\bigl(1-k^2\sn^2(t,k)\bigr),\\
\bigl({\textstyle\frac{\d}{\d t}}\cn(t,k)\bigr)^2&=\bigl(1-\cn^2(t,k)\bigr)
\bigl(1-k^2+k^2\cn^2(t,k)\bigr),\\
\bigl({\textstyle\frac{\d}{\d t}}\dn(t,k)\bigr)^2&=-\bigl(1-\dn^2(t,k)\bigr)
\bigl(1-k^2-\dn^2(t,k)\bigr),
\end{align*}
with initial conditions
\begin{alignat*}{3}
\sn(0,k)&=0,\;\> & \cn(0,k)&=1,\;\> & \dn(0,k)&=1,\\
{\textstyle\frac{\d}{\d t}}\sn(0,k)&=1,\;\>&
{\textstyle\frac{\d}{\d t}}\cn(0,k)&=0,\;\>&
{\textstyle\frac{\d}{\d t}}\dn(0,k)&=0.
\end{alignat*}

They satisfy the identities
\begin{equation*}
\sn^2(t,k)+\cn^2(t,k)=1 \;\>\text{and}\;\> k^2\sn^2(t,k)+\dn^2(t,k)=1,
\end{equation*}
and the differential equations
\e
\begin{gathered}
{\textstyle\frac{\d}{\d t}}\sn(t,k)=\cn(t,k)\dn(t,k),\qquad
{\textstyle\frac{\d}{\d t}}\cn(t,k)=-\sn(t,k)\dn(t,k)\\
\text{and}\qquad {\textstyle\frac{\d}{\d t}}\dn(t,k)=-k^2\sn(t,k)\cn(t,k).
\end{gathered}
\label{ru7eq4}
\e
When $k=0$ or 1 they reduce to trigonometric functions:
\begin{alignat*}{2}
\sn(t,0)&=\sin t,&\quad \cn(t,0)&=\cos t,\quad \dn(t,0)=1,\\
\sn(t,1)&=\tanh t&,\quad \cn(t,1)&=\dn(t,1)=\sech t.
\end{alignat*}
For $k\in[0,1)$ the Jacobi elliptic functions $\sn(t,k),\cn(t,k)$ 
and $\dn(t,k)$ are {\it periodic} in $t$, with a common period.

Using this notation we have the following result, adapted 
from~\cite[Th.~8.7]{Joyc1}.

\begin{thm} Let\/ $b_1,b_2,b_3$ be coprime integers satisfying $b_2\ge 
b_3>0>b_1$ and\/ $b_1+b_2+b_3=0$. Define $a>0$ and\/ $b\in[0,1)$ by 
\e
a^2=b_2(b_3-b_1)
\quad\text{and}\quad 
b^2=\frac{b_1(b_2-b_3)}{b_2(b_1-b_3)}.
\label{ru7eq5}
\e
Define $\phi:\R^2\ra{\cal S}^5$ by
\e
\begin{split}
\phi:(s,t)\mapsto\Bigl(i{\rm e}^{ib_1s}
\bigl({\ts\frac{b_2}{b_2-b_1}}\bigr)^{1/2}\dn(at,b),
&i{\rm e}^{ib_2s}\bigl({\ts\frac{b_1}{b_1-b_2}}\bigr)^{1/2}\cn(at,b),\\
&i{\rm e}^{ib_3s}\bigl({\ts\frac{b_1}{b_1-b_3}}\bigr)^{1/2}\sn(at,b)\Bigr),
\end{split}
\label{ru7eq6}
\e
and let\/ $N_0=\bigl\{r\,\phi(s,t):r,s,t\in\R\bigr\}$. Then $N_0$
is a special Lagrangian cone in $\C^3$ and\/ $\phi$ is an oriented
conformal map, so that\/ $(s,t)$ are oriented conformal coordinates
on $N_0\cap{\cal S}^5$. Furthermore, $\phi$ is doubly periodic in
$\R^2$, so that\/ $N_0$ is a two-sided\/ $T^2$-cone.
\label{ru7thm3}
\end{thm}

Applying Theorem \ref{ru6thm1} to this example with $w=u\frac{\pd}{\pd s}
+v\frac{\pd}{\pd t}$ for $u,v\in\R$, and using \eq{ru7eq4} to calculate 
$\psi=u\frac{\pd\phi}{\pd s}+v\frac{\pd\phi}{\pd t}$, yields:

\begin{thm} Let\/ $b_1,b_2,b_3$ be coprime integers satisfying $b_2\ge 
b_3>0>b_1$ and\/ $b_1+b_2+b_3=0$. Define $a>0$ and\/ $b\in[0,1)$ by 
\eq{ru7eq5}. Let\/ $u,v\in\R$ and define $N_{u,v}$ to be
\begin{align*}
\Bigl\{\Bigl(&\bigl({\ts\frac{b_2}{b_2-b_1}}\bigr)^{1/2}{\rm e}^{ib_1s}
\bigl((ir-ub_1)\dn(at,b)-ivab^2\,\sn(at,b)\cn(at,b)\bigr),\\
&\bigl({\ts\frac{b_1}{b_1-b_2}}\bigr)^{1/2}{\rm e}^{ib_2s}
\bigl((ir-ub_2)\cn(at,b)-iva\,\sn(at,b)\dn(at,b)\bigr),\\
&\bigl({\ts\frac{b_1}{b_1-b_3}}\bigr)^{1/2}{\rm e}^{ib_3s}
\bigl((ir-ub_3)\sn(at,b)+iva\,\cn(at,b)\dn(at,b)\bigr)\Bigr):
r,s,t\in\R\Bigr\}.
\end{align*}
Then $N_{u,v}$ is a ruled special Lagrangian $3$-fold in~$\C^3$.
\label{ru7thm4}
\end{thm}

It can be shown that when $b_1$ is even the 3-folds $N_{u,v}$
result from Theorem \ref{ru6thm3}, and are generically nonsingular
and diffeomorphic to $T^2\t\R$ as immersed submanifolds, and when
$b_1$ is odd and $u=0$ the 3-folds $N_{0,v}$ result from Theorem 
\ref{ru6thm4}, and are generically nonsingular and diffeomorphic 
to the total space of a real line bundle over the Klein bottle as 
immersed submanifolds.

We gave the cones of Theorem \ref{ru7thm3} as examples because
one can write them down in a very explicit way. But these are only
the simplest cases of two much larger explicit families of special 
Lagrangian $T^2$-cones in $\C^3$, which were constructed by the 
author in \cite[\S 8]{Joyc1} and \cite[\S 6]{Joyc2}, and which 
intersect in the examples of Theorem~\ref{ru7thm3}. 

The first of these families was also studied by Haskins 
\cite[\S 3--\S 5]{hask2}, and in terms of minimal Lagrangian tori in 
$\CP^2$ by Castro and Urbano \cite{CaUr}, and the second is related 
to examples due to Lawlor and Harvey, and was also studied by Bryant 
\cite[\S 3.5]{Brya2} from a different point of view. By applying Theorem 
\ref{ru6thm1} to these families we can obtain many more explicit examples 
of ruled SL 3-folds diffeomorphic to $T^2\t\R$ or a real line bundle over 
the Klein bottle. 

The constructions of \cite[\S 6]{Joyc2} also included a 1-parameter 
family of ruled SL 3-folds asymptotic to each $T^2$-cone, and results 
analogous to Theorems \ref{ru6thm3} and \ref{ru6thm4} are given for 
them in \cite[Th.s 6.3 \& 6.4]{Joyc2}. This 1-parameter family is part 
of the 2-parameter family of Theorem \ref{ru6thm3}, those with~$u=0$.

Finally, we briefly discuss $\U(1)$-invariant ruled SL 3-folds in $\C^3$.
Let $G\subset\SU(3)$ be a Lie subgroup isomorphic to $\U(1)$. What can
we say about $G$-invariant ruled SL 3-folds? Calculations by the author
show the following. There is a 2-dimensional family of $G$-invariant
SL cones, written down explicitly by Haskins \cite[\S 3--\S 5]{hask2} 
and the author \cite[\S 8]{Joyc1}. Applying Theorem \ref{ru6thm1}, we 
can enlarge this to a 4-dimensional family of explicit $G$-invariant
ruled SL 3-folds.

However, the family of all $G$-invariant ruled SL 3-folds is 6-dimensional. 
In the notation of \S\ref{ru5}, $\phi$ is already explicitly known by work 
of Haskins and the author, and we seek $G$-invariant solutions to the linear 
equations on $\psi$. These can be reduced to a linear first-order o.d.e.\ in 
4 variables. 

The coefficients of this o.d.e.\ involve the Jacobi elliptic functions, 
as these enter the explicit form of $\phi$. Two solutions to this o.d.e.\ 
are known from Theorem \ref{ru6thm1}, but the author has not been able 
to find the other two solutions, and so find an explicit form for general
$G$-invariant ruled SL 3-folds.

\end{document}